\documentclass[12pt]{amsart}
\usepackage{amssymb}
\usepackage{lscape}
\usepackage[headings]{fullpage}
\usepackage{amsfonts}
\usepackage{paralist}
\usepackage{xspace}
\usepackage{euscript}
\usepackage{color}
\usepackage{epigraph}
\usepackage[all]{xy}
\usepackage[normalem]{ulem}
\usepackage[colorlinks,pagebackref=true,pdftex]{hyperref}

\makeatletter
\@addtoreset{equation}{section}
\def\theequation{\thesection.\@arabic \c@equation}

\def\theenumi{\@roman\c@enumi}
\def\@citecolor{blue}
\def\@linkcolor{blue}
\def\@urlcolor{blue}

\makeatother

\newtheorem{lemma}[equation]{Lemma}
\newtheorem{prop}[equation]{Proposition}
\newtheorem{cor}[equation]{Corollary}

\newtheorem{claim*}{Claim}
\newtheorem{thm}[equation]{Theorem}

\theoremstyle{definition}
\newtheorem{remark}[equation]{Remark}
\newenvironment{rmk}[1][]{\begin{remark}[#1] \pushQED{\qed}}{\popQED \end{remark}}
\newtheorem{eg}[equation]{Example}
\newenvironment{example}[1][]{\begin{eg}[#1] \pushQED{\qed}}{\popQED \end{eg}}
\newtheorem{definition}[equation]{Definition}
\newenvironment{defn}[1][]{\begin{definition}[#1]\pushQED{\qed}}{\popQED \end{definition}}
\newtheorem{notn}[equation]{Notation}
\newenvironment{notation}[1][]{\begin{notn}[#1]\pushQED{\qed}}{\popQED \end{notn}}

\def\<{\langle}
\def\>{\rangle}

\newcommand{\codim}{\operatorname{codim}}

\newcommand{\coker}{\operatorname{coker}}

\newcommand{\Hom}{\operatorname{Hom}} 
\newcommand{\sheafHom}{\mathcal{H}om}

\newcommand{\Proj}{\operatorname{Proj}}

\newcommand{\rank}{\operatorname{rank}}

\newcommand{\Spec}{\operatorname{Spec}}

\newcommand{\Tor}{\operatorname{Tor}}

\renewcommand{\to}{\longrightarrow}

\newcommand{\bA}{\mathbb{A}}
\newcommand{\bb}{\mathbf{b}}
\newcommand{\bB}{\mathbf{B}}

\newcommand{\CC}{\mathbb{C}}
\newcommand{\EE}{\mathcal{E}}

\newcommand{\bg}{\mathbf{g}}
\newcommand{\KK}{\mathcal{K}}
\newcommand{\LL}{\mathcal{L}}
\newcommand{\fm}{\mathfrak m}

\newcommand{\NN}{\mathbb{N}}
\newcommand{\cO}{{\mathcal{O}}}
\newcommand{\PP}{\mathbb{P}}
\newcommand{\QQ}{\mathbb{Q}}

\newcommand{\bS}{\mathbf{S}}
\newcommand{\cS}{\mathcal{S}}
\newcommand{\Sc}{\mathrm{S}}
\newcommand{\rH}{\mathrm{H}}
\newcommand{\ES}{\mathrm{ES}}
\newcommand{\Di}{\mathrm{D}}
\newcommand{\wD}{\widetilde{\mathrm{D}}}
\newcommand{\cT}{\mathcal{T}}

\newcommand{\ZZ}{\mathbb{Z}}
\newcommand{\PPb}{\mathbb{P}(\vec{B})}
\newcommand{\GL}{{\bf GL}}

\DeclareMathOperator{\Seg}{Seg}

\newcommand{\defi}[1]{{\bfseries\upshape #1}}
\newcommand{\defeq}{:=}

\newcommand{\Bracket}[1]{\mbox{$\left[\begin{matrix}  #1\end{matrix}
\right]$}}
\newcommand{\minus}{\ensuremath{\!\smallsetminus\!}}
\newcommand{\pdim}{\operatorname{pdim}}
\newcommand{\EN}{\operatorname{EN}}

\newcommand{\excise}[1]{}

\newcommand{\arxiv}[1]{\href{http://arxiv.org/abs/#1}{{\tt arXiv:#1}}}

\title[Tensor complexes]%
{Tensor complexes: %
Multilinear free resolutions constructed from higher tensors}

\author[C. Berkesch Zamaere]{Christine Berkesch Zamaere}
\address{Institut Mittag-Leffler \\ Aurav\"agen 17 \\ SE-182 60 Djursholm, Sweden \hfill \quad \linebreak\vspace{-.7cm}}
\address{
Department of Mathematics \\ Stockholm University \\
SE-106 91 Stockholm, Sweden}
\curraddr{School of Mathematics \\ University of Minnesota \\ Minneapolis, MN 55455} 
\email{cberkesc@math.umn.edu}

\author[D. Erman]{Daniel Erman}
\address{Department of Mathematics \\ Stanford University \\
Stanford, CA 94305}
\curraddr{Department of Mathematics \\ University of Wisconsin \\ Madison, WI 53706}
\email{derman@math.wisc.edu}

\author[M. Kummini]{Manoj Kummini}
\address{Department of Mathematics \\ Purdue University \\
West Lafayette, IN 47907}
\curraddr{Chennai Mathematical Institute \\ Siruseri, Tamilnadu, 603103. India.}
\email{mkummini@cmi.ac.in}

\author[S. Sam]{Steven V Sam}
\address{Department of Mathematics \\ Massachusetts Institute of
Technology \\ Cambridge, MA 02139}
\curraddr{Department of Mathematics \\ University of California \\ Berkeley, CA 94720}
\email{svs@math.berkeley.edu}

\thanks{The first author was supported by NSF Grant 
OISE 0964985.  The second author was partially supported 
by an NDSEG fellowship and NSF Award No. 1003997.  
The fourth author was supported by an NSF graduate research fellowship and an NDSEG fellowship.}


\begin{document}
\begin{abstract}
  The most fundamental complexes of free modules over a
  commutative ring are the Koszul complex, which is constructed from a vector (i.e., a
  $1$-tensor), and the Eagon--Northcott and Buchsbaum--Rim
  complexes, which are constructed from a matrix (i.e.,  a $2$-tensor).  
  The subject of this paper is a multilinear
  analogue of these complexes, which we construct from an arbitrary
  higher tensor.  

  Our construction provides detailed new examples of minimal free
  resolutions, as well as a unifying view on a wide variety of
  complexes including: the Eagon--Northcott, Buchsbaum--Rim and
  similar complexes, the Eisenbud--Schreyer pure resolutions, and the
  complexes used by Gelfand--Kapranov--Zelevinsky and Weyman to
  compute hyperdeterminants.  In addition, we provide applications to
  the study of pure resolutions and Boij--S\"oderberg theory, 
  including the construction of infinitely many new families of pure
  resolutions, and the first explicit description of the differentials
  of the Eisenbud--Schreyer pure resolutions.
  \end{abstract}

\maketitle

\vspace{-.3cm}
\section{Introduction}
\label{sec:intro}
 \begin{epigraphs}
 \qitem{In commutative algebra, the Koszul complex is the mother of all complexes.}%
       {David Eisenbud}
 \end{epigraphs}
 The most fundamental complex of free modules over a commutative ring
 $R$ is the Koszul complex, which is constructed from a vector (i.e.,
 a $1$-tensor) $\mathbf{f}=(f_1, \dots, f_a)\in R^a$.  The next most
 fundamental complexes are likely the Eagon--Northcott and
 Buchsbaum--Rim complexes, which are constructed from a matrix (i.e.,
 a $2$-tensor) $\widetilde{\psi}\in R^a\otimes R^b$.

In this paper we construct 
multilinear analogues of these complexes, which we refer to as
\defi{tensor complexes}. These complexes are constructed
from an arbitrary
higher tensor $\widetilde{\phi}\in R^a\otimes R^{b_1} \otimes
\dots\otimes R^{b_n}$, 
providing a unifying perspective on
many of these previously known families --- including Koszul,
Eagon--Northcott, and Buchsbaum--Rim complexes --- 
and leads to new
such families of resolutions.  This also supplies a new tool
 for producing and studying invariants of higher tensors.

While tensor complexes display remarkable
numerical properties (for instance, all extremal rays of the cone of
Betti diagrams can be generated by our construction; see
\S\ref{sec:ESpures}), their structure is surprisingly
simple. We provide explicit descriptions of
these free resolutions from several different perspectives; in
particular, each tensor complex can be pieced together from
linear strands of a Koszul complex.
This not only adds tensor complexes to the few families 
of free resolutions that are understood in detail, 
it also provides new such families 
that are uniformly minimal over $\ZZ$.  (Uniformity over $\ZZ$ can
be quite subtle; see~\cite{hashimoto}.)

To motivate our main result, 
we first recall some properties of the more familiar Eagon--Northcott
complex. 
The Eagon--Northcott complex for an arbitrary matrix can be
constructed as a pullback from the universal case. 
Namely, if we first build the Eagon--Northcott complex $\EN(\psi)_\bullet$ 
over the polynomial ring $\ZZ[x_{i,j}]$ for an $a\times b$ matrix 
$\psi=\psi^{a\times b}=(x_{i,j})$ of indeterminates,
then the Eagon--Northcott complex of
$\widetilde{\psi}$ is $\EN(\psi)_\bullet \otimes_{\ZZ[x_{i,j}]} R$.
Several nice properties of the complex $\EN(\psi)_\bullet$ 
are illustrated in the following theorem.

\begin{thm}[{Eagon--Northcott \cite{eagon-northcott}}]
\label{thm:ENcomplex}
The Eagon--Northcott complex $\EN(\psi)_\bullet$ of a matrix of indeterminates
$\psi$ satisfies the following:
\begin{compactenum}[\rm (i)]
\item  It is a graded free resolution of a Cohen--Macaulay module.
\item It is uniformly minimal over $\ZZ$,
	i.e., $\EN(\psi)_\bullet \otimes_{\ZZ[x_{i,j}]} \Bbbk[x_{i,j}]$ is a
	minimal free resolution for any field $\Bbbk$. 
\item It is a pure resolution, 
	i.e., $\EN(\psi)_i$ is generated in a single degree for each $i$.
\item It respects the bilinearity of $\psi$, 
	i.e., $\EN(\psi)_\bullet$
	is $\GL_a\times \GL_b$-equivariant.
\end{compactenum}
\end{thm}

The Buchsbaum--Rim complex also satisfies the assertions of
Theorem~\ref{thm:ENcomplex}.  In fact, the Eagon--Northcott and
Buchsbaum--Rim complexes fit naturally into a sequence of bilinear
complexes arising from the matrix $\psi$ and a weight $w\in
\ZZ^2$~\cite{buchs-eis}.\footnote{\cite[\S A2.6]{eisenbud} outlines
  the construction of matrix complexes, and we use this as our primary
  reference for these complexes.  There, the complexes are
  parametrized by $\ZZ^1$, which corresponds to the second coordinate
  of our $w\in\ZZ^2$; the first coordinate of $w$ simply allows a
  twist of the complex as a whole.}  We refer to an element of this
sequence as a \defi{matrix complex}, although these are sometimes
called ``generalized Koszul complexes'' (see~\cite{buchsbaum,
  buchsbaum-rim}). While such a complex exists for any $w$, an
analogue of Theorem~\ref{thm:ENcomplex} holds only for a limited set
of weights.

To construct the tensor complexes of an arbitrary tensor
$\widetilde{\phi}$, we similarly take the pullback of the
universal case.  Let $a\in \NN$ and $\bb=(b_1, \dots, b_n)\in \NN^n$.
We define a universal tensor $\phi:=\phi^{a\times \bb}$ over the
symmetric algebra $S=\Sc^{\bullet}(\ZZ^a\otimes \ZZ^{b_1}\otimes \dots
\otimes \ZZ^{b_n})$, and in \S\ref{subsec:Fphiw:defn}, we construct
the \defi{tensor complex} $F(\phi,w)_\bullet$ from this universal 
tensor and a weight $w\in \ZZ^{n+1}$. 
The following theorem illustrates how tensor complexes are 
a multilinear extension of the Eagon--Northcott complex and 
the other matrix complexes, as long as we limit the
choice of $w$, requiring it to be a \defi{pinching weight} (see
Definition~\ref{def:pinching:weight}).

\begin{thm}\label{thm:tensor:complexes}
If $w$ is a pinching weight for $\phi^{a\times \bb}$, then 
$F(\phi,w)_\bullet$ satisfies the following:
\begin{compactenum}[\rm (i)]
\item \label{item:tensor:CM} It is a graded free resolution of
  a Cohen--Macaulay module $M(\phi,w)$. 
\item \label{item:tensor:flat} 
	It is uniformly minimal over $\ZZ$.
\item\label{item:tensor:pure} 
	It is a pure resolution.
\item \label{item:tensor:eqvt} 
	It respects the multilinearity of $\phi$, 
	i.e., $F(\phi,w)_\bullet$ is 
	$\GL_a \times \dots \times \GL_{b_n}$-equivariant.
\end{compactenum}
\end{thm}

A connection between tensors and free complexes has
previously been observed in special cases,~\cite[\S 14]{gkz} and
\cite[\S 9.4]{weyman}. For instance,~\cite[Proposition~14.3.2]{gkz}
uses a free complex to express hyperdeterminants of the boundary format, 
and this is a special case of our
construction (see Proposition~\ref{prop:hyperdeterminant}).
Hyperdeterminants also play an important role 
in the study of general tensor complexes.  
As shown in Theorem~\ref{thm:hyperdeterminants2}, 
the support of $M(\phi,w)$ is set-theoretically defined
by an ideal of hyperdeterminants of certain sub-tensors of $\phi$. 
In addition, each such variety is a resultant variety for a system of multilinear equations on a product of projective spaces 
(see Proposition~\ref{prop:resultant}).

Tensor complexes extend another important
class of free resolutions: pure resolutions of Cohen--Macaulay modules.  
Such resolutions are central objects in Boij--S\"oderberg theory, 
as they provide the extremal rays of the cone of Betti diagrams. 
We show in Theorem~\ref{thm:family}
that there are an infinite number of different
tensor complexes whose Betti diagrams lie on any such
extremal ray.  In addition, Theorem~\ref{thm:ESpures} shows that each
Eisenbud--Schreyer pure resolution
from ~\cite[\S 5]{ES-JAMS} is obtained by taking hyperplane sections 
of a tensor complex, thus providing the first explicit description of the 
differentials of these complexes.

Properties of higher tensors are the subject of much recent work
(see~\cite{landsberg-complexity,landsberg-book} for surveys).  
A tensor complex for an arbitrary tensor $\widetilde{\phi}$ 
attaches new invariants to the tensor. In some small cases (see
Example~\ref{ex:rank}), these invariants detect the
rank of the tensor.  It would be interesting to pursue further connections.

\subsection{Constructing tensor complexes}
\label{subsec:tensorCxConst}
Perhaps the most important feature about the tensor complexes 
$F(\phi,w)_\bullet$ is that we can describe them explicitly. 
To underscore their essential properties, we 
present three different perspectives on these complexes.

\subsubsection*{Strands of the Koszul complex}
\label{subsec:strands}
In~\cite[\S A2.6]{eisenbud}, matrix complexes are constructed by splicing
together two strands of a Koszul complex. Tensor complexes are similar: if $\phi$ is an $(n+1)$-tensor
and $w$ is a pinching weight for $\phi$, then $F(\phi,w)_\bullet$ can be
built by splicing together $n$ strands of a Koszul complex. 

For example, consider the universal $7\times (2,2)$ tensor $\phi=\phi^{7\times (2,2)}$.
Let $A\cong \ZZ^7$, $B_1\cong \ZZ^2\cong B_2$, $X^{7\times (2,
  2)}:=A\otimes B_1^*\otimes B_2^*$, and $S:=\Sc^{\bullet}(X^{7\times
  (2,2)})$.
For the choice of pinching weight $w=(0,1,4)$, the tensor complex
$F(\phi,w)_\bullet$ is 
\[
\xymatrix{
S^{10} &  S^{28}(-1) \ar[l] &  S^{70}(-3) \ar[l]_-{\sigma} &  
S^{70}(-4)\ar[l] &  S^{28}(-6)\ar[l]_-{\sigma'} &  S^{10}(-7)\ar[l] & 0\ar[l]
}.
\]
To illustrate the equivariant structure of this free resolution,
in \S\ref{subsec:notation:rep:theory} we introduce a column notation 
for writing representations of $\GL(A)\times \GL(B_1)\times \GL(B_2)$, 
giving $F(\phi,w)_\bullet$ the form 
\begin{equation}\label{eqn:threestrand}
\begin{split}
{
\footnotesize
\xymatrix @C=7mm @R=2mm{
\Bracket{\wedge^0\\ \Sc^1\\ \Sc^4}
&\Bracket{\wedge^1\\ \Sc^0\\ \Sc^3}(-1)\ar[l]
\\
&&
\Bracket{{\wedge}^{3}\\ \wD^0\\ \Sc^1}(-3)\ar[ul]_-{\sigma}
&\Bracket{{\wedge}^{4}\\ \wD^1 \\ \Sc^0}(-4)\ar[l]
\\
&&&& 
\Bracket{{\wedge}^{6}\\ \wD^3 \\ \wD^0}(-6)\ar[ul]_-{\sigma'}
&\Bracket{{\wedge}^{7}\\ \wD^4\\ \wD^1}(-7)\ar[l]
&0.\ar[l]
}}
\end{split}
\end{equation}
Here, for instance, the $F_2$ term of \eqref{eqn:threestrand} 
denotes the graded free $S$-module
$\bigwedge^3(A)\otimes_{\ZZ} \det(B_1^*)\otimes_{\ZZ} \Sc^1(B_2)\otimes_{\ZZ} S(-3)$. 
This complex arises from three separate strands --- vertically 
separated in \eqref{eqn:threestrand} --- of a Koszul
complex $K(\phi)_\bullet$ on the $\ZZ^3$-graded polynomial ring
$\Sc^{\bullet}(X^{7\times (2,2)}\otimes B_1\otimes B_2)$.  While this
mirrors the construction of matrix complexes in \cite[\S A2.6]{eisenbud},
it will be modified in our situation by certain
local cohomology modules (see \S\ref{sec:strands}). 
We splice these strands together via the maps $\sigma$ and $\sigma'$ 
whose entries are expressions in the $2\times 2$ minors of the flattening 
$\phi^{\flat} \colon A^*\otimes S \to B^*_1\otimes B^*_2\otimes S$.
The fact that $F(\phi,w)_\bullet$ forms a complex then follows from a
generalized Laplace expansion formula for the determinant of a singular
matrix.  Example~\ref{ex:422} provides a detailed illustration of this
fact in a similar example.

For a tensor complex, a new phenomenon arises that was not present in the case of matrix complexes: it is possible that two consecutive
maps are splice maps.  In fact, there will be many cases where none of
the differentials $F(\phi,w)_\bullet$ consist of linear forms; 
each strand consists of a single free module and each 
differential is a splicing map.

\subsubsection*{Tensor complexes and representation theory}
\label{subsec:reptheory}
The above approach to $F(\phi,w)_\bullet$ makes little use of
its multilinear symmetry.  By incorporating ideas from representation 
theory, we are able to provide a simple description of the differentials of
$F(\phi,w)_\bullet$.

Let us reconsider the map $\sigma$ from \eqref{eqn:threestrand}.  This map
is determined by its degree 3 part $[\sigma]_3\colon [F_2]_3\to
[F_1]_3$, which is the following map of finite-rank free $\ZZ$-modules:
\[
[\sigma]_3 \colon
\wedge^3 A\otimes \wedge^2B_1^*\otimes \Di^0B_1^*\otimes \Sc^1B_2 
\to  
\left( \wedge^1 A \otimes \Sc^0B_1\otimes \Sc^3B_2 \right)\otimes
\Sc^2(X^{7\times (2,2)}).
\]
Recalling that $X^{7\times (2,2)}=A\otimes B_1^*\otimes B_2^*$, we
express the map $[\sigma]_3$ entirely in terms of tensor products and adjoints of multiplication and comultiplication maps.
Namely, we use the subrepresentation $\wedge^2 A \otimes \wedge^2
B_1^*\otimes \Di^2 B_2^* \subseteq \Sc^2(X^{7\times (2,2)})$ 
and construct $[\sigma]_3$ via the following equivariant maps
on each tensor factor:
\[
[\sigma]_3 \leftrightarrow
\begin{cases}
  \wedge^3A \to \wedge^1A \otimes \wedge^2A & \text{by
  comultiplication,}\\
  \wedge^2B_1^* \otimes \Di^0B_1^* \to \Sc^0B_1\otimes
  \wedge^2B_1^* & \text{by identifying } \Di^0B_1^*\cong \Sc^0B_1,\\
  \Sc^1B_2\to \Sc^3B_2\otimes \Di^2B_2^* & \text{by the adjoint
    of multiplication}.
\end{cases}
\]
This provides an explicit description of the differentials of
$F(\phi,w)_\bullet$ (see \S\ref{sec:differentials}) and proves that
$F(\phi,w)_\bullet$ is a complex (see Lemma~\ref{lemma:directcomplex}).
For acyclicity, we take a third perspective.

\subsubsection*{The geometric method}
\label{subsec:geomethod}
The geometric method of
Kempf--Lascoux--Weyman~\cite[\S 5]{weyman} provides the most powerful perspective for studying the tensor complex
$F(\phi,w)_\bullet$.  Continuing with the
universal $7\times (2,2)$ tensor example, we define a complex
$\KK(\phi,w)_\bullet$ on $\Spec(S)\times \PP(B_1)\times \PP(B_2)$ as
the sheafy version of $K(\phi)_\bullet$, twisted by a line bundle
determined by $w$.  Taking the derived pushforward of
$\KK(\phi,w)_\bullet$ along the projection $\pi \colon
\Spec(S)\times \PP(B_1)\times\PP(B_2)\to \mathbb \Spec(S)$
also yields the tensor complex $F(\phi,w)_\bullet$; 
we use this as our primary definition of $F(\phi,w)_\bullet$ (see
Definition~\ref{defn:Fphiw}).

The geometric method immediately
provides the acyclicity of $F(\phi,w)_\bullet$. The disadvantage is that
the geometric method does not provide a clear description
of the differentials of the complex. To make use
of the representation theoretic description in \S\ref{sec:differentials},
it suffices to show that the differentials can be chosen equivariantly.
(This is not obvious, since the
representation theory of $\GL_n(\ZZ)$ is not semisimple.)

\subsection{The algebra and geometry of tensor complexes}
\label{subsec:further}
We now summarize some additional results on tensor complexes, 
as well as applications of our work to Boij--S\"oderberg theory. 
We begin with the functorial properties of $F(\phi,w)_\bullet$. 

\begin{prop}\label{prop:functor}
  Let $a' \leq a$, and let $w$ and $w'$ be weights. Let $S :=
  \ZZ[X^{a\times\bb}]$ and $S' := \ZZ[X^{a'\times \bb}]$.  
Given an inclusion $i \colon \ZZ^{a'} \to \ZZ^a$ and a polynomial of
  multi-degree $w-w'$ in $S'\otimes \Sc^\bullet (B_1)\otimes \dots
  \otimes \Sc^\bullet(B_n)$, we have a degree zero map of complexes
  \[
  F(\phi^{a'\times \bb},w')_\bullet\otimes_{S'} S \to F(\phi^{a\times
    \bb},w)_\bullet.
  \]
\end{prop}

This result is proven in \S\ref{sec:functoriality} and
is related to~\cite[Theorem~1.2]{BEKS}, as the maps considered
in that result are special cases of the above construction.

We now turn to properties of the module $M(\phi,w)$ 
that is resolved by the tensor complex $F(\phi,w)_\bullet$;
these statements are proved in \S\ref{sec:Mphi}.

\begin{cor} \label{cor:Mphiw} Let $\phi=\phi^{a\times \mathbf{b}}$ be
  the universal tensor and $w$ be a pinching weight for $\phi$.
\begin{enumerate}[\rm (i)]
	\item  \label{enum:SupportMphiw}
    The support of $M(\phi,w)$ is an irreducible subvariety of $\mathbb
    A^{a\times \bb}$ that is independent of $w$ and has 
    codimension $a-\sum_{i} (b_i-1)$. 
    \item  \label{enum:MphiwGenPerf}
    $M(\phi,w)$ is
    generically perfect, i.e., it is Cohen--Macaulay and faithfully
    flat over $\ZZ$. 
	\item  \label{enum:MphiwMult} 
   The multiplicity of $M(\phi,w)$ is 
   independent of $w$. Specifically, it is given by the multinomial coefficient
   \[
   e( M(\phi,w))= \frac{a!}{(a -\sum_i (b_i-1))!  \prod_{i=1}^n
     (b_i-1)!}.
   \]
\end{enumerate}
\end{cor}

\subsubsection*{Hyperdeterminantal varieties}
\label{subsec:supports}
Based on Corollary~\ref{cor:Mphiw}\eqref{enum:SupportMphiw}, 
we denote the support of $M(\phi,w)$ by $Y(\phi)$ 
and call such a variety a \defi{hyperdeterminantal variety}. 
These hyperdeterminantal varieties simultaneously extend 
the determinantal varieties defined by maximal minors of a matrix 
(when $\phi$ is a $2$-tensor)
and the hypersurfaces defined by hyperdeterminants of the boundary
format (see~\cite[\S 14.3]{gkz}). 

\begin{thm}\label{thm:hyperdeterminants2}
Let $\phi=\phi^{a\times \mathbf{b}}$ and $Y(\phi)\subseteq \mathbb A^{a\times \bb}$ be the support variety of $M(\phi,w)$. 
If $a':=1+\sum_{i=1}^n (b_i-1)$, then $Y(\phi)$ is
set-theoretically defined by the ideal 
\begin{equation}\label{eq:hyperdet ideal}
\left\< \text{\emph{hyperdeterminant of }} \phi'  \mid  \phi' \text{ is an } (a'\times \bb)-\text{subtensor of  } \phi \right\>.
\end{equation}
\end{thm}
The ideal~\eqref{eq:hyperdet ideal} can fail to be radical, 
as we illustrate in Example~\ref{ex:422support}.
Further, Remark~\ref{rmk:resultant} explains how 
the variety $Y(\phi)$ is a resultant variety for a system of multilinear equations on a product of projective spaces, yielding the following result.

\begin{prop}\label{prop:resultant}
For a field $\Bbbk$, let $\mathbf{f}=f_1, \dots, f_a$ 
be a collection of multilinear forms on 
$\PP_\Bbbk^{b_1-1}\times \dots \times \PP_\Bbbk^{b_n-1}$. 
This gives a tensor 
$\phi_{\mathbf{f}}
\in \Bbbk^a\otimes \Bbbk^{b_1}\otimes \dots \otimes \Bbbk^{b_n}$ 
and thus a specialization map $q_{\mathbf{f}} \colon S\to \Bbbk$, 
which sends $\phi\mapsto \phi_{\mathbf{f}}$. 
Let $w$ be any pinching weight for the universal tensor $\phi^{a\times \bb}$, 
and let $\partial_\bullet$ denote the differential of $F(\phi^{a\times\bb},w)$. 
Denote by $\partial_1(\mathbf{f})$ the matrix obtained by specializing 
the entries of $\partial_1$ via the map $q_{\mathbf{f}}$.  
The following are then equivalent:
\begin{enumerate}[\rm (i)]
	\item The vanishing locus 
	$V(f_1,\dots,f_a)\subseteq 
	\PP^{b_1}_\Bbbk\times \dots \times \PP^{b_n-1}_\Bbbk$ is nonempty 
    \textup{(}over any algebraic closure of $\Bbbk$\textup{)}.
	\item  The matrix $\partial_1(\mathbf{f})$ does not have full rank.
\end{enumerate}
\end{prop}
We explore the geometry of hyperdeterminantal varieties
in \S\ref{sec:MDVs}.  In contrast to the case of determinantal
varieties, we show the varieties $Y(\phi)$ are rarely normal or
Cohen--Macaulay.

\subsubsection*{Applications to Boij--S\"oderberg theory}
\label{subsec:appsBStheory}
The construction of tensor complexes has significant implications for
Boij--S\"oderberg theory (see~\cite{ES:ICMsurvey} for a survey) and the study of pure resolutions.  
A sequence $d=(d_0,\dots,d_p)\in \ZZ^{p+1}$ is a \defi{degree
  sequence} if $d_i<d_{i+1}$ for all $i$. For a degree sequence $d$,
we say that $G_\bullet$ is a \defi{pure resolution of type} $d$ if for
each $i$, $G_i$ is generated in degree $d_i$.

\begin{thm}\label{thm:family}
  Let $d=(d_0, \dots, d_p)\in \ZZ^{p+1}$ be a degree sequence.  Then
  there exist infinitely many choices of $a, \bb$, and $w$ such that
  $w$ is a pinching weight for $\phi^{a \times \bb}$, $F(\phi^{a\times
    \bb},w)_\bullet$ is a pure resolution of type $d$, and
  $M(\phi^{a\times \bb},w)$ is a Cohen--Macaulay module that is flat
  over $\ZZ$.
\end{thm}

The pure resolutions of type $d$ constructed
in Theorem~\ref{thm:family} are unrelated to one another; this 
yields infinitely many new families of pure resolutions of type $d$ for
every $d$.
More precisely,
we have the following:
\begin{prop}\label{prop:indecomp}
Suppose that $a \geq \sum_{j=1}^n (b_j-1)$ and that $w$ is a pinching
weight for $\phi^{a\times \bb}$. Then $M(\phi^{a\times \bb},w)$ is
indecomposable.
\end{prop}

Theorem~\ref{thm:family} builds on previous work of
\cite{efw,ES-JAMS}.  
Namely, two constructions of Cohen--Macaulay
modules with a pure resolution of type $d$ were previously known: in
characteristic 0 given in \cite[\S\S 3,4]{efw} and a different
construction that works in arbitrary characteristic in 
\cite[\S 5]{ES-JAMS}.  The Eisenbud--Schreyer construction 
arises as a hyperplane section of a certain tensor
complex (see Theorem~\ref{thm:ESpures}).  However, this is 
unsurprising, as our original motivation for this project was to
understand a multilinear version of their work.

Our results thus provide the first explicit description of pure
resolutions over a field of positive characteristic, as we produce 
a closed formula for their differentials 
(without the need to explicitly compute the pushforward of a complex).  
In characteristic zero, a similar explicit description for the pure
resolutions of \cite[\S 3]{efw} appears
in~\cite[\S\S1,2]{pieri-resolutions}.  In another direction, a
recent algorithm of Eisenbud, based on the
Bernstein--Gelfand--Gelfand correspondence, 
enables the computation of the differentials of the pushforward of 
a complex. This algorithm would compute the differentials of
any specific pure resolution of~\cite[\S 5]{ES-JAMS}
and is implemented in~\cite[\texttt{BGG} package, version~1.4]{M2}.


Finally, we note that the construction of the tensor 
complex $F(\phi,w)_\bullet$ extends
  to any scheme.  Namely, if $\widetilde{\phi}$ is a global section of
  a tensor product of vector bundles $\mathcal A\otimes \mathcal
  B_1\otimes \cdots \otimes \mathcal B_n$ on a scheme $X$, then there
  is a natural $\cO_X$-module version of the complex
  $F(\phi,w)_\bullet$.

\subsection{Outline}
\label{subsec:outline}
We outline our notation in \S\ref{sec:notation} and describe the
general geometric construction of the complex $F(\phi,w)_\bullet$. In
\S\ref{sec:balanced}, we introduce a particularly nice class of tensor
complexes, called ``balanced tensor complexes'' and discuss their
basic properties. The differentials of these complexes are described
explicitly using representation-theoretic methods in 
\S\ref{sec:differentials}.

Beginning with \S\ref{sec:main}, we turn our attention to the main
construction of tensor complexes, proving
Theorem~\ref{thm:tensor:complexes} via
Theorem~\ref{thm:main:structure}.  \S\ref{sec:strands} describes the
construction of tensor complexes from strands of the Koszul complex,
and \S\ref{sec:functoriality} illustrates the functorial properties of
tensor complexes.  \S\S\ref{sec:Mphi} and \ref{sec:MDVs} examine
properties of the modules $M(\phi,w)$ and their supports $Y(\phi)$,
respectively.

In \S\ref{sec:ESpures}, we relate the Eisenbud--Schreyer construction
of pure modules to balanced tensor complexes.  Further
applications of our main results to Boij--S\"oderberg theory,
including the construction of new families of pure resolutions, can be
found in \S\ref{sec:new:pures}.  

Finally, \S\ref{sec:detailed} provides a detailed example of a 
tensor complex, including presentation matrices for the differentials.
We have also provided Appendix~\ref{sec:charfree}, which reviews some
basic definitions and constructions in multilinear algebra, and
Appendix~\ref{sec:schurfunctors} provides a rapid review of the
facts we employ from the representation theory 
of the general linear group over a field of characteristic zero.

\subsection*{Acknowledgements}
We thank D.~Eisenbud and J.~Weyman for many thoughtful 
discussions about this work.
We also thank W.~Heinzer, D.~Katz, S.~Kleiman, B.~Sturmfels, 
and T.~V\'{a}rilly-Alvarado for helpful comments.  
This work began during the ``Workshop on Local Rings and Local
Study of Algebraic Varieties'' at ICTP,
was continued during the AMS Mathematics Research Community on 
``Commutative Algebra'' and a workshop funded 
by the Stanford Mathematics Research Center, 
and was completed while the first author attended the program 
``Algebraic Geometry with a view towards applications" at Institut Mittag-Leffler;
we are grateful for all of these opportunities.
Throughout the course of this work, calculations were performed 
using the software {\tt Macaulay2} \cite{M2}.

\section{Notation and general construction of the complex 
$F(\phi,w)_\bullet$}
\label{sec:notation}
Let $a\in \NN$ and $\bb=(b_1, \dots, b_n)\in\NN^n$. 
Let $A$ and $B_j$, $1 \leq j \leq n$, be free $\ZZ$-modules of
rank $a$ and $b_j$, $1 \leq j \leq n$, respectively.  
We define $\bB:=B_1\otimes_{\ZZ} \dots \otimes_{\ZZ} B_n$ and
$X := X^{a\times \bb}:= A\otimes \bB^*$.  

For a free $\ZZ$-module $V$, we write $\ZZ[V]$ for the symmetric
algebra on $V$. Throughout, $S$ will denote the polynomial ring
$S:=\ZZ[X^{a\times\bb}]=\ZZ[x_{i,J}]$, where $1\leq i\leq a$ and $J=(j_1, \dots, j_n)$, 
$1\leq j_\ell \leq b_\ell$.  We endow $S$ with the standard
$\ZZ$-grading $\deg(x_{i,J})=1$.  We write the universal tensor
$\phi=\phi^{a\times \bb}$ as
\[
\phi=(x_{i,J})\in S\otimes_{\ZZ} X^{a\times \bb}.
\]

Given an $(n+1)$-tensor, there are a number of ways to obtain a matrix
by ``flattening'' this tensor.  One such flattening is particularly
useful for our purposes.  Via the isomorphism $S\otimes_{\ZZ} X^{a\times
  \bb}\cong \Hom_S(S\otimes_{\ZZ}A^*,S\otimes_{\ZZ} \bB^*)$, $\phi$ induces
  a map of free $S$-modules
\[
\phi^{\flat} \colon S \otimes_\ZZ A^* \to S \otimes_\ZZ \bB^*.
\]

We write $\PP(B_j)$ for the projective space of $1$-dimensional
quotients of $B_j$, so that $\PP(B_j)=\Proj (\ZZ[B_j])$.  Let $\PPb :=
\PP(B_1) \times \cdots \times \PP(B_n)$ and 
$\bA^{a\times \bb} \defeq \Spec( \ZZ[X^{a\times \bb}])$.  

\subsection{Representation theory conventions}
\label{subsec:notation:rep:theory}
Let $G = \GL(A) \times \GL(B_1) \times \cdots
\times \GL(B_n)$.  For a free $\ZZ$-module $V$ of finite rank, we use
$\Sc^i(V)$ to refer to its $i$th symmetric power, 
$\Di^i(V)$ for its $i$th divided power, 
and $\det(V)$ for its top exterior power.  
We are most interested in divided powers 
twisted by a copy of the determinant, so we set
\[
\wD^i(V):=\Di^i(V)\otimes \det(V).
\]
We use the convention that $\rH^0(\PP^0_{\ZZ},\cO(d))\cong \Sc^d(\ZZ)$
for all $d$. Although $\GL_1(\ZZ) \cong \ZZ/2$ cannot distinguish
between two different $d$ of the same parity, these representations
are distinct from a ``functor of points'' perspective, i.e., they are
distinct over larger coefficient rings, such as $\QQ$. Similar remarks
apply to powers of the determinant representation in general. 
When $V$ is a $\QQ$-vector space, we use ${\bf S}_\lambda V$ 
to denote irreducible representations of $\GL(V)$.  See Appendix~\ref{sec:schurfunctors}
for a summary of representation theory results used in this paper.

We write the representations over $G$ as columns, so that the order of
the rows allows us to omit the reference to the free modules $A, B_1,
\dots, B_n$. Inside the columns, we abbreviate $\wD^i(B^*)$ as
$\wD^i$.  A twist by $S(-i)$ is denoted by $(-i)$ next to the column.
For example,
\[
\Bracket{\wedge^0\\ \Sc^1\\ \Sc^4} :=
\wedge^0(A)\otimes \Sc^1(B_1)\otimes \Sc^4(B_2) \otimes S
\quad \text{ and } \quad 
\Bracket{\wedge^3\\ \wD^1\\ \Sc^1}(-2) := 
\wedge^3(A)\otimes \wD^1(B_1^*)\otimes \Sc^1(B_2) \otimes S(-2).
\]

\subsection{Free resolution conventions}
\label{subsec:resolutions}
Conventions for the graded Betti diagrams of graded free complexes are
standard. Namely, let $L_\bullet$ be a graded free complex over
$S$. The graded Betti numbers $\beta_{i,j}(L_\bullet)$ are defined as
follows:
\[
L_i=\bigoplus_{j\in \ZZ} S(-j)^{\beta_{i,j}(L_\bullet)}.
\]
The Betti diagram of $L_\bullet$ is then 
\[
\beta(L_\bullet)=\begin{pmatrix}
\vdots & \vdots &&\vdots \\
\beta_{0,-1} & \beta_{1,0} & \cdots & \beta_{p,p-1}\\
\beta_{0,0} & \beta_{1,1} & \cdots  & \beta_{p,p}\\
\beta_{0,1} & \beta_{1,2} &\cdots & \vdots \\
\vdots &\vdots & \ddots &
\end{pmatrix}.
\]
Betti diagrams have nonzero entries in only finitely many positions, 
so we omit the rows of zeroes in examples.

\begin{defn}\label{defn:unifMinl}
  Let $M$ be a finitely generated graded $S$-module and 
$L_\bullet$ a free resolution of $M$.
  We say that $L_\bullet$ is \defi{uniformly minimal} if
  $L_\bullet \otimes_\ZZ \Bbbk$ is a minimal free resolution for every
  field $\Bbbk$. In this case, we define
  $\beta_{i,j}(M):=\beta_{i,j}(L_\bullet)$.
\end{defn} 

\subsection{General construction of $F(\phi,w)_\bullet$}
\label{subsec:Fphiw:defn}
To construct the complex $F(\phi,w)_\bullet$ we apply a minor
extension of the geometric method~\cite[\S 5]{weyman}, working
over $\ZZ$ instead of an arbitrary field.  For the reader unfamiliar with~\cite{weyman},
this may be a rather opaque definition.  
Several concrete descriptions of these complexes are given later (see
Proposition~\ref{prop:balanced}, Definition~\ref{defn:differentials},
Theorem~\ref{thm:main:structure}, and \S\ref{sec:strands}).

To apply this extension of the geometric method, we observe that 
the lemmas in \cite[\S 5.2]{weyman} hold over $\ZZ$ 
if the sheaves involved are flat over $\ZZ$ 
and all of the relevant sheaf cohomology 
(as in \eqref{eqn:cohom:in:each:deg}) is free over $\ZZ$. 
In our situation, this is the case. 
(Alternatively, one can prove acyclicity of the relevant complexes 
over $\ZZ$ by proving acyclicity over each finite field as well as $\QQ$, 
in which case the results of \cite[\S 5]{weyman} apply directly.)

Recall that $\PPb = \PP(B_1) \times \cdots \times \PP(B_n)$, and view 
$\bA^{a\times \bb}\times \PPb$ as the total space of the trivial
bundle $\EE:=A^*\otimes_{\ZZ} \bB\otimes_{\ZZ} \cO_{\PPb}$ over
$\PPb$.
Consider the vector bundle $\cT := A^* \otimes \cO_{\PPb}(1, 1,\dots, 1)$ on $\PPb$.
There is a natural surjective map $\EE\to \cT$ induced by
the natural maps $B_i\otimes \cO_{\PP(B_i)}\to
\cO_{\PP(B_i)}(1).$ Let $\cS$ be the kernel of this map, so
that we have an exact sequence of vector bundles on $\PPb$ of the form 
\[
\xymatrix{%
  0\ar[r] &  \cS \ar[r] & \EE \ar@2{-}[d] \ar[r] & \cT\ar[r]
  \ar@2{-}[d]&  0 \\ 
  && A^*\otimes_\ZZ \bB\otimes_\ZZ \cO_{\PPb} & A^*\otimes_\ZZ
  \cO_{\PPb}(1,\dots, 1) }.
\]
Explicitly, $\cS = \sheafHom((\bB^* \otimes_\ZZ \cO_{\PPb}) /
\cO_{\PPb}(-1, \dots, -1), A^* \otimes_\ZZ \cO_{\PPb})$, and 
we let
$Z(\phi) = Z(\phi^{a\times\bb}) \subseteq \bA^{a\times \bb} \times \PPb$
denote $\mathbb V_{\PPb}(\cS)$, the total space  of $\cS$. 
The total space $\mathbb V_{\PPb}(\EE)$ of the vector bundle $\EE$ 
is $\bA^{a\times \bb}\times \PPb$. 
Write $\pi \colon \bA^{a\times \bb}\times \PPb\to \bA^{a\times \bb}$ 
for the projection. Let
$Y(\phi) = \pi(Z(\phi))$ scheme-theoretically. Note that $Z(\phi)$
and $Y(\phi)$ are integral schemes. We have a commutative diagram:
\begin{equation} \label{eqn:Yphi}
  \xymatrix{
    Z(\phi)=\mathbb V_{\PPb}(\cS) \ar[rr] \ar[d]^{\mu} & &
    \mathbb A^{a\times \bb}\times \PPb=\mathbb V_{\PPb}(\EE)\ar[d]^{\pi}\\
    Y(\phi)\ar[rr]&&\mathbb A^{a\times \bb}.
  }
\end{equation}
Let $\pi_2 \colon \bA^{a\times \bb}\times \PPb\to \PPb$ be the natural
projection, and consider the following Koszul complex on
$\bA^{a\times\bb}\times \PPb$,
\begin{equation}\label{eq:sheaf:Koszul}
\KK(\phi)_\bullet:  \quad
 \cO_{ \bA^{a\times \bb}\times \PPb}
\longleftarrow
 \bigwedge^1(\pi_2^*\cT^*)
\longleftarrow
 \cdots 
\longleftarrow
\bigwedge^{a-1}(\pi_2^*\cT^*)
\longleftarrow
 \bigwedge^{a}(\pi_2^* \cT^*) 
\longleftarrow
0
\end{equation}
which resolves the sheaf $\cO_{Z(\phi)}$. 

\begin{defn} \label{defn:Fphiw}
Fix a \defi{weight vector} $w=(w_0,\dots, w_n)\in\ZZ^{n+1}$. To define $F(\phi,w)_\bullet$, we follow the construction of ~\cite[Theorem~5.1.2]{weyman}, with $\cO_{\PPb}(w_1,\dots, w_n)$ in place of $\mathcal V$. Additionally, we twist
the resulting complex by $S(-w_0)$ to obtain the \defi{tensor complex}, denoted $F(\phi,w)_\bullet$. 

It follows immediately that $F(\phi,w)_\bullet$ is a graded, free complex 
of $S$-modules that is quasi-isomorphic to 
${\bf R}\pi_* \left( \KK(\phi)_\bullet\otimes 
\pi_2^*\cO_{\PPb}(w_1, \dots, w_n)\right) \otimes_S S(-w_0)$.
The terms of $F(\phi,w)_\bullet$ are
\begin{align}\label{eqn:cohom:in:each:deg}
  F(\phi,w)_i & = \bigoplus_{j \ge 0} \rH^j(\PPb, \bigwedge^{i+j} \cT^*
  \otimes \cO_{\PPb}(w_1, \dots, w_n))
  \otimes S(-i-j-w_0) \nonumber \\
  & = \bigoplus_{j \ge 0} \rH^j(\PPb, \cO_{\PPb}(w_1-i-j, \dots, w_n-i-j))
  \otimes \bigwedge^{i+j} A \otimes S(-i-j-w_0).
\end{align}

We write $\partial_i$ for the differential $F(\phi,w)_i \to F(\phi, w)_{i-1}$. 
Let $M(\phi,w) := \coker \partial_1$.
\end{defn}
There is a minor abuse of notation inherent in the above definition.
Namely, to define the differentials of such a complex via the
geometric method, we must explicitly compute a free complex to
represent the quasi-isomorphism class of a pushforward of
a complex, and there is some choice involved in building this complex
(see~\cite[\S 5.5]{weyman}).
Thus the differentials $\partial_i$ are not a priori
determined by $\phi$ and $w$.  We ignore this subtlety because 
our main cases of interest are when $w$ is a
pinching weight for $\phi$, and in these cases, we may make a
canonical choice for each differential (up to sign) via representation
theory, as illustrated in Proposition~\ref{prop:canonical} and
Theorem~\ref{thm:main:structure}.

\begin{rmk}\label{rmk:Koszul:forms}
  Let $K(\phi)_\bullet$ denote the $\ZZ^{n+1}$-graded complex of
  graded free $\ZZ$-modules on
  $\ZZ[X^{a\times\bb}]\otimes\ZZ[B_1]\otimes \dots \otimes \ZZ[B_n] = \ZZ[x_{i,J},y_{j_\ell}]$
  corresponding to the Koszul complex of sheaves $\KK(\phi)_\bullet$
  from \eqref{eq:sheaf:Koszul}.  For $J=(j_1, \dots, j_n)$, set $y_J:=y_{j_1}\cdots y_{j_n}$.  Consider the multilinear forms
  \[
  f_i:=\sum_J x_{i,J}y_J, \quad i = 1, \ldots, a.
  \]
  Then $K(\phi)_\bullet$ is the
  Koszul complex on $(f_1, \dots, f_a)$.
\end{rmk}

\begin{rmk} \label{rmk:super} 
If we replace $A$ by any $\ZZ/2$-graded free $\ZZ$-module and 
take care in using $\ZZ/2$-graded multilinear algebra 
(see, for example, \cite[\S 2.4]{weyman}), 
essentially all of our assertions about tensor complexes remain true, 
with one significant difference. 
If the odd part of $A$ is nonzero, then $S$ will be a graded 
commutative algebra, and the resulting complexes will be 
infinite in length in one direction. If the even part of $A$ is 0, 
then we obtain pure resolutions over the exterior algebra. 
\end{rmk}

\section{Balanced tensor complexes}
\label{sec:balanced} 

In \S\ref{sec:notation} we defined $F(\phi,w)_\bullet$ for an
arbitrary weight vector $w\in\ZZ^{n+1}$.  To obtain free resolutions
with nice properties, including those
outlined in Theorem~\ref{thm:tensor:complexes}, we impose further
conditions on the weight vector $w$.  For clarity, we begin by 
introducing a particularly simple class of examples called 
balanced tensor complexes.
The construction is sufficiently rich to produce tensor complexes that 
are pure resolutions of type $d$ for every degree sequence $d$. 
In fact, this construction is closely modeled on the 
Eisenbud--Schreyer construction of pure resolutions 
\cite[\S 5]{ES-JAMS}.  In \S\ref{sec:main}, we extend the
results of this section to more general tensor complexes.

\begin{defn}\label{defn:balanced}
We say that $F(\phi^{a\times \bb},w)_\bullet$ is a 
\defi{balanced tensor complex} if it satisfies the following conditions:
\begin{enumerate}
\item \label{item:balanced:sum} $a=b_1+b_2+\cdots +b_n$.
\item \label{item:balanced:weights} $w_1=0$ and $w_i=b_1+\dots +b_{i-1}$ for
  $i=2,\dots, n$.
\end{enumerate}
Set $d(w):=(w_0, w_2 + w_0, w_3 + w_0, \dots, w_{n-1} + w_0, w_n +
w_0, a + w_0)\in \ZZ^{n+1}$.
\end{defn}

The condition \eqref{item:balanced:sum} is less restrictive than 
it appears because we allow 
the possibility of tensoring with rank-$1$ free modules.
For instance, there is a natural way to
identify a $7\times (3,2)$ tensor with a $7\times (3,1,2,1)$ tensor or
with a $7\times (1,1,3,2)$ tensor and so on. 
These identifications enable us to produce many examples 
of balanced tensor complexes. 
The following example illustrates this flexibility. 

\begin{example}[Complexes of~{\cite[\S A2.6]{eisenbud}}] \label{ex:matrix:complexes} Let $b \leq a \in \NN$.
The matrix complexes $\EuScript{C}^0, \ldots,
  \EuScript{C}^{a-b}$ of~\cite[\S A2.6]{eisenbud} may be realized
  as examples of balanced tensor complexes. Fix $0 \leq i
  \leq a-b$.  Let 
\[
\bb := 
(\underbrace{1, \ldots, 1}_{i}, b, \underbrace{1, \ldots, 1}_{a-b-i}).
\]
The corresponding balanced tensor complex $F(\phi^{a\times \bb},
  w)_\bullet$ is isomorphic to $\EuScript{C}^{i}$.
\end{example}

The following proposition proves a portion of 
Theorem~\ref{thm:tensor:complexes} for balanced tensor complexes.

\begin{prop}\label{prop:balanced}
Suppose that 
$F(\phi,w)_\bullet$ is a balanced tensor complex. 
Write $d(w) = (d_0, \ldots, d_n)$. Then
\begin{equation}\label{eqn:balancedMods}
F(\phi,w)_i\cong S(-d_i)\otimes \bigwedge^{d_i-w_0} A \otimes
\bigotimes_{j=1}^i \wD^{d_{i}-d_{j}} (B_j^*) \otimes
\bigotimes_{j=i+1}^n \Sc^{d_{j-1}-d_{i}}(B_j). 
\end{equation}
In particular, $F(\phi,w)_\bullet$ is a pure resolution of type $d(w)$ and
satisfies Theorem~\ref{thm:tensor:complexes}(\ref{item:tensor:CM}--\ref{item:tensor:pure}).
\end{prop}

\begin{proof}
From \eqref{eqn:cohom:in:each:deg}, we must consider sheaves of 
the form
$\bigwedge^l \cT^* \otimes \cO(w_1,\dots,w_n) 
= \bigwedge^l A \otimes \cO(w_1-l,\dots,w_n-l)$, 
which is nonzero only if $l \in [0,a]$.
By the K\"unneth formula, this sheaf will have nonzero
cohomology precisely when $l\notin [w_i+1,w_i+b_i-1]$ for all $i =
1,\dots,n$.  Since $w_{i+1}=w_i+b_i$, it immediately follows that $l\in
\{0=w_1,w_2,\dots,w_n,w_n+b_n=a\}$.

Set $d':=(0,w_2,\dots,w_n,w_n+b_n)$, and note that $d_i'+w_0=d_i$ for all $i$. 
Computing the cohomology for $l=d_i'$ yields
\begin{align*}
    F(\phi,w)_i &= S(-d_i) \otimes \bigwedge^{d_i-w_0} A \otimes
    \bigotimes_{j=1}^i \rH^{b_j - 1}(\PP(B_j), \cO(w_j - d'_i)) \otimes 
    \bigotimes_{j=i+1}^n \rH^0(\PP(B_j), \cO(w_j - d'_i)),
  \end{align*}
  which is~\eqref{eqn:balancedMods}.  In particular, the complex has
  no terms in negative homological degrees, and hence
  \cite[Theorem 5.1.2]{weyman} implies that $F(\phi,w)_\bullet$ is a
  minimal free resolution of $M(\phi,w)$ and $M(\phi,w)\otimes_S
  S(w_0)$ is naturally isomorphic to $\rH^0(\PPb, \Sc^\bullet(\cS^*) \otimes_\ZZ
  \cO_{\PPb}(w_1,\dots,w_n))$.  Since the latter is free over $\ZZ$,
  $M(\phi,w)$ is also free over $\ZZ$, completing the proof of 
  Theorem~\ref{thm:tensor:complexes}\eqref{item:tensor:flat}.

  We now prove that $M(\phi,w)$ is Cohen--Macaulay (i.e.,
  Theorem~\ref{thm:tensor:complexes}\eqref{item:tensor:CM}). 
  Since we know that $\pdim M(\phi,w) = n \geq \codim M(\phi,w)$, 
  it suffices to show that $\codim M(\phi,w) \geq
  n$.  By~\cite[Theorem~5.1.2(b)]{weyman}, the support of $M(\phi,w)$
  is the variety $Y(\phi)$ from \eqref{eqn:Yphi}.  Recall that
  $Z(\phi)$ is the total space of $\cS$.  The codimension of $Z(\phi)$
  in $\mathbb A^{a\times \bb}\times \PPb$ thus equals the rank of 
  $\cT$, which is $a$. Therefore 
  \[
  \dim Y(\phi) \le \dim Z(\phi)= 
  \dim \bA^{a \times \bb} + \dim \PPb - a=\dim \bA^{a \times \bb}-n,
  \]
  so $\codim Y(\phi) \ge n$, as desired.
\end{proof}

We provide a more detailed description of the support of $M(\phi, w)$
in~\S\ref{sec:Mphi}.  We also note that for any degree sequence $d$ 
there exists a unique balanced tensor complex $F(\phi,w)_\bullet$ 
that is a pure resolution of type $d$. 
This follows from Theorem~\ref{thm:ESpures}.

\begin{example} \label{ex:034711} 
Take $a=11$ and $\bb = (3,1,3,4)$. To obtain a balanced
complex, we set $w = (0,0,3,4,7)$. Then $d(w)=(0,3,4,7,11)$ with 
the following free resolution:
\[
\xymatrix{
S^{1800} &  S^{17325}(-3) \ar[l] &  S^{19800}(-4) \ar[l] &  
S^{4950}(-7)\ar[l] &  S^{675}(-11)\ar[l]& 0\ar[l]
},
\]
which we denote:
\[
\small \xymatrix @-0.5pc {
\Bracket{\wedge^0 \\ \Sc^0 \\ \Sc^3 \\ \Sc^4 \\ \Sc^7} 
&
\Bracket{\wedge^3 \\ \wD^0 \\ \Sc^0 \\ \Sc^1 \\ \Sc^4}(-3)\ar[l]
&
\Bracket{\wedge^4 \\ \wD^1 \\ \wD^0 \\ \Sc^0 \\ \Sc^3}(-4)\ar[l]
&
\Bracket{\wedge^7 \\ \wD^4 \\ \wD^3 \\ \wD^0 \\ \Sc^0}(-7)\ar[l]
&
\Bracket{\wedge^{11} \\ \wD^8 \\ \wD^7 \\ \wD^4 \\ \wD^0}(-11)\ar[l]
& 
0 \ar[l]
}. \qedhere
\]
\end{example}

\section{Explicit differentials for balanced tensor complexes}
\label{sec:differentials}

Since our definition of $F(\phi,w)_\bullet$ involves an application of
the geometric method, we would a priori need to
explicitly compute the pushforward of a complex in order to define a specific
differential.  In this section, we use a representation theoretic
argument to illustrate that such a computation is unnecessary
for balanced tensor complexes. 
Definition~\ref{defn:differentials} describes the equivariant
differential, and the main result of
this section is Proposition~\ref{prop:canonical}.  In
\S\ref{sec:main}, we extend this proposition to the more general
setting of Theorem~\ref{thm:tensor:complexes}.

\begin{prop}\label{prop:canonical}
Let $F(\phi,w)_\bullet$ be a balanced tensor complex.  Up to sign, there is
a unique differential $\partial_\bullet$, defined explicitly in
Definition~\ref{defn:differentials}, which makes $F(\phi,w)_\bullet$ into a
$G$-equivariant free resolution.  In particular,
$(F(\phi,w)_\bullet,\partial_\bullet)$ satisfies 
Theorem~\ref{thm:tensor:complexes}\eqref{item:tensor:eqvt}.
\end{prop}

For a free module $F_i$ we use $[F_i]_{e}$ to denote the degree $e$
piece of $F_i$; for a map of free modules $f \colon F_i\to F_{i-1}$ we
use $[f]_e$ to denote the induced map $[f]_e\colon [F_i]_e\to
[F_{i-1}]_e$.  To define the $G$-equivariant differentials
$\partial_i\colon F(\phi,w)_i \to F(\phi,w)_{i-1}$, we define a
$G$-equivariant map $[\partial_i]_{d_i} \colon [F(\phi,w)_i]_{d_i} \to
[F(\phi,w)_{i-1}]_{d_i}$ on the generators of $F(\phi,w)_i$ and
extend $S$-linearly.

By Proposition~\ref{prop:balanced}, the source of $[\partial_i]_{d_i}$
is given by
\begin{align} \label{eq:decomp:domain:splice} 
  [F(\phi,w)_i]_{d_i} = \bigwedge^{d_i - w_0} A \otimes
  \bigotimes_{j=1}^i \wD^{d_i - d_j}(B_j^*) \otimes
  \bigotimes_{j=i+1}^n \Sc^{d_{j-1} - d_i}(B_j).
\end{align}
Noting that $b_i = d_i - d_{i-1}$, the corresponding decomposition for
the target is
\begin{equation} \label{eq:decomp:image:splice}
  \begin{split}
    [F(\phi,w)_{i-1}]_{d_i} 
    &= [F(\phi,w)_{i-1}]_{d_{i-1}} \otimes \Sc^{b_i}(X)\\
    &= \left( \bigwedge^{d_{i-1} - w_0} A \otimes
      \bigotimes_{j=1}^{i-1} \wD^{d_{i-1} - d_j}(B_j^*) \otimes
      \bigotimes_{j=i}^n \Sc^{d_{j-1} - d_{i-1}}(B_j) \right) \otimes
    \Sc^{b_i}(X).
  \end{split}
\end{equation}

By Appendix~\ref{sec:charfree}, we have inclusions of $G$-modules:
\begin{align} 
\label{eqn:subreprOfSym}
\begin{split}
  \Sc^{b_i}(X) 
  &\supseteq \bigwedge^{b_i} A \otimes
  \bigwedge^{b_i}(B_1^* \otimes \cdots \otimes B_n^*) \\
  &\supseteq \bigwedge^{b_i} A \otimes \det B_i^* \otimes
  \Di^{b_i}(\bigotimes_{j \ne i} B_j^*)\\
  &\supseteq \bigwedge^{b_i} A \otimes \det B_i^* \otimes 
  	\bigotimes_{j\ne i} \Di^{b_i}(B_j^*). 
\end{split}
\end{align}

\begin{defn}[Equivariant Differentials on
  $F(\phi,w)_\bullet$]\label{defn:differentials}
  The map $[\partial_i]_{d_i} \colon [F(\phi,w)_i]_{d_i} \to
  [F(\phi,w)_{i-1}]_{d_i}$ is defined to be the composition of a map
  \[
  \iota\colon [F(\phi,w)_i]_{d_i} \longrightarrow
  [F(\phi,w)_i]_{d_{i-1}}\otimes \left( \bigwedge^{b_i} A \otimes \det
    B_i^* \otimes \bigotimes_{j\ne i} \Di^{b_i}(B_j^*)\right)
  \]
  with the inclusion obtained from \eqref{eqn:subreprOfSym}. We define
  $\iota$ to be the tensor product $\iota=\iota_A\otimes
  \iota_{B_1}\otimes \dots \otimes \iota_{B_n}$ where the components
  are defined below.  For $\iota_A$ we take the comultiplication
  map
  \[
  \iota_A \colon \bigwedge^{d_i-w_0} A \longrightarrow
  \bigwedge^{d_{i-1}-w_0} A \otimes \bigwedge^{b_i} A.
  \]
  For $j\leq i-1$, we take the twist by $\det(B_j^*)$ of the dual of
  the multiplication map $\Sc^{b_i}(B_j)\otimes
  \Sc^{d_{i-1}-d_j}(B_j)\to \Sc^{d_i-d_j}(B_j)$, and set:
  \[
  \iota_{B_j} \colon \wD^{d_i - d_j} (B_j^*) \longrightarrow
  \wD^{d_{i-1} - d_j}(B_j^*) \otimes \Di^{b_i}(B_j^*).
  \]
  For $j=i$, we choose an identification (unique up to sign)
  $\wD^0(B_i^*)\cong \Sc^0(B_i)\otimes \det(B_i^*)$. Finally, when
  $j\geq i+1$ we take the dual of the contraction map $\Di^{d_{j-1} -
    d_{i-1}}(B_j^*)\otimes \Sc^{b_i}(B_j)\to \Di^{d_{j-1}-d_i}(B_j^*)$, and
  set:
 \[
 \iota_{B_j} \colon \Sc^{d_{j-1} - d_i}(B_j) \longrightarrow
	\Sc^{d_{j-1} - d_{i-1}}(B_j) \otimes \Di^{b_i}(B_j^*).
 \]
We then define $\partial_i \colon F(\phi,w)_i\to F(\phi,w)_{i-1}$ as the
$S$-linear extension of $[\partial_i]_{d_i}$.  
The map $\partial_i$ is clearly $G$-equivariant.
\end{defn}

We say that a map of free $\ZZ$-modules is \defi{saturated} if its
cokernel is also a free $\ZZ$-module.

\begin{lemma}\label{lem:saturated}
The map $[\partial_i]_{d_i}$ is saturated and injective.
\end{lemma}
\begin{proof}
  Since $[\partial_i]_{d_i}$ is the tensor product of the maps
  $\iota_A$ and $\iota_{B_j}$, it suffices to show that each of these
  maps is saturated and injective.  For $\iota_A$ this follows
  from~\cite[Theorems III.1.4, III.2.4]{abw}. For $j\ne i$, the map
  $\iota_{B_j}$ is the dual of a surjective map of free $\ZZ$-modules,
  so it is saturated and injective. Finally, the map $\iota_{B_i}$ is an
  isomorphism.
\end{proof}

The following lemma is essential to the claim of uniqueness in
Proposition~\ref{prop:canonical}.
\begin{lemma}[Base change to $\QQ$]\label{lemma:multfree}
  The $G(\QQ)$-representation
  \begin{align*}
  [F(\phi,w)_i]_{d_i}\otimes \QQ
  &= \bigwedge^{d_i - w_0} A \otimes
  \bigotimes_{j=1}^i \wD^{d_i - d_j}(B_j^*) \otimes
  \bigotimes_{j=i+1}^n \Sc^{d_{j-1} - d_i}(B_j) \otimes \QQ
\end{align*}
appears with multiplicity 1 inside $[F(\phi,w)_{i-1}]_{d_i} \otimes
\Sc^{b_i}(X) \otimes \QQ$.
\end{lemma}  
\begin{proof}
  We first find the subrepresentations $W = \bS_\lambda A \otimes
  \bS_{\mu^1} B_1^* \otimes \cdots \otimes \bS_{\mu^n} B_n^*\otimes
  \QQ$ of $\Sc^{b_i}(X)\otimes \QQ$ whose tensor product with
  $[F(\phi,w)_{i-1}]_{d_{i-1}}\otimes \QQ$ contains
  $[F(\phi,w)_i]_{d_i}\otimes \QQ$. By Pieri's rule \eqref{eqn:pieri},
  this only happens for $\lambda = \mu^i = (1^{b_i})$ and $\mu^j =
  (b_i)$ for $j \ne i$.  By Schur--Weyl duality
  (\eqref{eqn:schurweyl1} and \eqref{eqn:schurweyl2}), $W$ appears in
  $\Sc^{b_i}(X)\otimes \QQ$ with multiplicity 1.
\end{proof}

There is a straightforward proof that $\partial^2=0$, which we
include below.

\begin{lemma}\label{lemma:directcomplex}
  For any $i\geq 1$ we have $\partial_{i}\partial_{i+1}=0$.  In
  particular, $(F(\phi,w)_\bullet, \partial_\bullet)$ is a complex.
\end{lemma}
\begin{proof}
  It is enough to verify that the composition
  $[F(\phi,w)_{i+2}]_{d_{i+2}} \to [F(\phi,w)_{i+1}]_{d_{i+2}} \to
  [F(\phi,w)_i]_{d_{i+2}}$ is 0 where $0\le i \le n-2$.  Since
  $F(\phi, w)_\bullet$ is a free complex, we may tensor with $\QQ$ before
  checking that this map is 0. Thus, for the rest of the proof, we
  assume that all free $\ZZ$-modules have been tensored by $\QQ$.

  Since the maps $\partial_i$ are $G$-equivariant, it suffices to show
  that any $G(\QQ)$-equivariant map $[F(\phi,w)_{i+2}]_{d_{i+2}} \to
  [F(\phi,w)_i]_{d_{i+2}}$ is zero.  First write $e := d_{i+2} - d_i =
  b_{i+2}+b_{i+1}$. By \eqref{eqn:cauchy},
  \[
  \Sc^eX = \bigoplus_{\lambda \vdash e} \bS_\lambda A \otimes
  \bS_\lambda(\bB^*).
  \]
  By Pieri's rule \eqref{eqn:pieri}, the generators of
  $[F(\phi,w)_{i+2}]_{d_{i+2}}$ can only appear in the tensor product
  of the generators of $[F(\phi,w)_i]_{d_{i+2}}$ with the direct
  summand with $\bS_\lambda = \bigwedge^e$.  Now by
  \eqref{eqn:cauchy},
  \[
  \bigwedge^e(\bB^*) = \bigoplus_{\mu \vdash e} \bS_\mu(B_1^*
  \otimes \cdots \otimes B_{i}^* \otimes B_{i+3}^* \otimes \cdots
  \otimes B_n^*) \otimes \bS_{\mu'}(B_{i+1}^* \otimes B_{i+2}^*).
  \]
  Using Pieri's rule \eqref{eqn:pieri} and Schur--Weyl duality
  (\eqref{eqn:schurweyl1} and \eqref{eqn:schurweyl2}), we only need to
  focus on the summand with $\bS_\mu = \Sc^e$. By \eqref{eqn:cauchy},
  \[
  \bigwedge^e (B_{i+1}^* \otimes B_{i+2}^*) = \bigoplus_{\nu \vdash e}
  {\bS}_\nu B_{i+1}^* \otimes \bS_{\nu'} B_{i+2}^*,
  \]
  so we must show that 
\[
(\bS_\nu B_{i+1}^* \otimes \bS_{\nu'} B_{i+2}^*) \otimes \Sc^{b_{i+1}} B_{i+2}
\]
  does not contain a copy of $(\Sc^{b_{i+2}} B_{i+1})^* \otimes \det
  B_{i+1}^* \otimes \det B_{i+2}^*$ for any $\nu \vdash e$.  Since $\rank
  B^*_{i+2} = b_{i+2}$, this happens precisely when $\nu'$ is the
  partition $(b_{i+1} + 1, 1^{b_{i+2} - 1})$.  However, in this case
  $\bS_\nu B_{i+1}^* = 0$ because $\rank B_{i+1} = b_{i+1}$.
\end{proof}

\begin{prop}
\label{prop:balanced:eqvt}
The complex $(F(\phi,w)_\bullet, \partial_\bullet)$ is a free
resolution of $M(\phi,w)$.
\end{prop}

\begin{proof}
  To simplify notation, we drop reference to $\phi$ and $w$ throughout
  this proof. Let $(F_\bullet, \epsilon_\bullet)$ be a uniformly
  minimal free resolution of $M$.  We use $\epsilon_0 \colon F_0\to M$
  to denote the natural quotient map. From
  Lemma~\ref{lemma:directcomplex}, $(F_\bullet, \partial_\bullet)$ is
  a free complex. We set $\partial_0:=\epsilon_0$.
  
  We first claim that $\partial_0\partial_1=0$. This can be checked
  after base changing to $\QQ$. By \cite[Theorem~5.4.1]{weyman}, the
  complex $F_\bullet \otimes \QQ$ admits a
  $G(\QQ)$-equivariant differential $\epsilon'_\bullet$ which makes it
  acyclic. By Lemma~\ref{lemma:multfree}, $[\partial_1]_{d_1}\otimes
  \QQ$ is a nonzero scalar multiple of $[\epsilon'_1]_{d_1}$, and thus
  $\epsilon'_0\epsilon'_1 = 0 = \partial_0\partial_1$.
 
  Now, since $(F_\bullet, \epsilon_\bullet)$ is a resolution of $M$ and
  $(F_\bullet, \partial_\bullet)$ is a free complex mapping to $M$ (by
  $\partial_0$), the identity $M\overset{\text{id}}\to M$ induces a map of
  complexes $a_\bullet\colon (F_\bullet, \partial_\bullet)\to (F_\bullet,
  \epsilon_\bullet)$ by~\cite[Lemma~20.3]{eisenbud}.  We claim that $a_i$
  is an isomorphism for each $i$, and we proceed by induction.

  For $i=0$, we may assume that $a_0$ is the identity. For the
  induction step, we assume that $a_i$ is an isomorphism, so we have a diagram:
  \[
  \xymatrix{ 0 \ar[r] & [F_{i+1}]_{d_{i+1}}
    \ar[rr]^-{[\epsilon_{i+1}]_{d_{i+1}}} && [F_i]_{d_{i+1}} \ar[rr] &&
    \coker\left([\epsilon_{i+1}]_{d_{i+1}} \right) \ar[r] &0\\
    0 \ar[r] & [F_{i+1}]_{d_{i+1}} \ar[rr]^-{[\partial_{i+1}]_{d_{i+1}}}
    \ar[u]_{[a_{i+1}]_{d_{i+1}}}&& [F_i]_{d_{i+1}}
    \ar[rr]\ar[u]_{[a_{i}]_{d_{i+1}}}^{\cong}&&
    \coker\left([\partial_{i+1}]_{d_{i+1}} \right)\ar[r] \ar[u]_{b}&0.
  }
  \] 
  Since the middle arrow is an isomorphism, it follows that $b$ is
  surjective.  The cokernel of $[\epsilon_{i+1}]_{d_{i+1}}$ is a free
  $\ZZ$-module since the complex $(F_\bullet, \epsilon_\bullet)$ is a
  uniformly minimal resolution, and the cokernel of
  $[\partial_{i+1}]_{d_{i+1}}$ is a free $\ZZ$-module by
  Lemma~\ref{lem:saturated}. Thus, $b$ is an isomorphism.  By the five
  lemma, we conclude that $[a_{i+1}]_{d_{i+1}}$ is an isomorphism of
  $\ZZ$-modules, and hence $a_{i+1}$ is an isomorphism of $S$-modules.
\end{proof}

We are now prepared to prove the main result of this section.

\begin{proof}[Proof of Proposition~\ref{prop:canonical}]
  Theorem~\ref{thm:tensor:complexes}\eqref{item:tensor:eqvt} follows from
  Proposition~\ref{prop:balanced:eqvt}.  For uniqueness, assume that
  $\epsilon_\bullet$ is another $G$-equivariant differential.
  Lemma~\ref{lemma:multfree}, after a base-change to $\QQ$, implies
  that $\epsilon_i$ and $\partial_i$ differ by an integer scalar
  multiple. By uniform minimality, this integer cannot be divisible by
  any prime number, so it must be $\pm 1$.
\end{proof}

\begin{rmk}[Kronecker coefficients]
In characteristic 0, the acyclicity of $(F(\phi, w)_\bullet, \partial_\bullet)$ 
imposes nonvanishing conditions on the Kronecker coefficients 
$g_{\lambda, \mu^1, \dots, \mu^n}$ 
(see Appendix~\ref{sec:schurfunctors} for the relevant definitions 
and results). 
For example, let $n=2$ and consider the first differential
  \[
  \bigwedge^{d_1} A \otimes \bigwedge^{d_1} B_1^* \otimes S(-d_1)
  \to \Sc^{d_1} B_2 \otimes S.
  \]
  When $i < d_2 - d_1$, the $G$-equivariant map
  \[
  \bigwedge^{d_1} A \otimes \bigwedge^{d_1} B_1^* \otimes \Sc^i(A
  \otimes B_1^* \otimes B_2^*) \to \Sc^{d_1} B_2 \otimes
  \Sc^{i+d_1}(A \otimes B_1^* \otimes B_2^*)
  \]
  is injective. Now rewrite the left-hand side as
  \[
  \bigwedge^{d_1} A \otimes \bigwedge^{d_1} B_1^* \otimes
  \bigoplus_{\lambda, \mu, \nu \vdash i} (\bS_\lambda A \otimes
  \bS_\mu B_1^* \otimes \bS_\nu B_2^*)^{\oplus g_{\lambda, \mu, \nu}},
  \]
  and the right-hand side as
  \[
  \Sc^{d_1} B_2 \otimes \bigoplus_{\alpha, \beta, \gamma \vdash
    i+d_1} (\bS_\alpha A \otimes \bS_\beta B_1^* \otimes \bS_\gamma
  B_2^*)^{\oplus g_{\alpha, \beta, \gamma}}.  
  \]
  It follows that if $g_{\lambda, \mu, \nu} \ne 0$, 
  then for any partition $\alpha$ obtained from $\lambda$ by adding a
  vertical strip of size $d_1$ and $\beta = (\mu_1 + 1, \dots,
  \mu_{d_1} + 1)$, there exists $\gamma$ obtained from $\nu$ 
  by adding a horizontal strip of size $d_1$ 
  such that $g_{\alpha, \beta, \gamma} \ne 0$.   
\end{rmk}

\subsection{Writing the differentials via minors of flattenings}
\label{subsec:flattenings}
Definition~\ref{defn:differentials} provides the following method for
writing the differentials of $F(\phi,w)_\bullet$ explicitly in terms of
minors of the flattening $\phi^{\flat}$.  If we choose bases $\{f_k\}$ and
$\{g_l\}$ of $F(\phi,w)_i$ and $F(\phi,w)_{i-1}$, then we may represent
$\partial_i$ by a matrix $\Psi$ of polynomials of degree $b_i$.
Consider the map 
\[
\alpha\colon [F(\phi,w)_i]_{d_i}\otimes
[F(\phi,w)_{i-1}]_{d_{i-1}}^*\longrightarrow \Sc^{b_i}(X),
\]
which is adjoint to $[\partial_i]_{d_i}$.
Note that $\alpha(f_k,g_l^*)$ is the $(k,l)$th entry of $\Psi$. 
Now consider the adjoint $\gamma$ of the map $\iota$ given in
Definition~\ref{defn:differentials}:
\[
\gamma \colon [F(\phi,w)_i]_{d_i} \otimes
[F(\phi,w)_{i-1}]_{d_{i-1}}^*\longrightarrow \bigwedge^{b_i} A \otimes
\det B_i^* \otimes \bigotimes_{j\ne i} \Di^{b_i}(B_j^*).
\]
Since $[\partial_i]_{d_i}$ was defined in terms of $\iota$ and the
inclusion~\eqref{eqn:subreprOfSym}, it follows that $\alpha$ is given by
$\gamma$ and~\eqref{eqn:subreprOfSym}.  

The first line of~\eqref{eqn:subreprOfSym} corresponds to the
inclusion of the $b_i\times b_i$ minors of $\phi^{\flat}$ into the
space of all polynomials of degree $b_i$.  Hence each entry of $\Psi$
may be defined in terms of $b_i\times b_i$ minors of $\phi^{\flat}$,
and we may write $\alpha$ explicitly via a formula for the inclusion
\[
 \bigwedge^{b_i} A \otimes \det B_i^* \otimes 
  	\bigotimes_{j\ne i} \Di^{b_i}(B_j^*) \subseteq 
	 \bigwedge^{b_i} A \otimes
  \bigwedge^{b_i}(B_1^* \otimes \cdots \otimes B_n^*).
\]
We obtain the necessary formula for this inclusion from repeated
applications of the multilinear inclusions described in
Appendix~\ref{sec:charfree}.

\begin{example}\label{ex:422inclusion}
  Let $a\times \bb=4\times (2,2)$ and $w=(0,0,2)$.  The complex
  $F(\phi,w)_\bullet$ has the form
\[
\xymatrix{
\Bracket{ \wedge^0 \\ \Sc^0\\ \Sc^2}&
\Bracket{\wedge^2 \\ \wD^0\\  \Sc^0}(-2)\ar[l]_-{\partial_1}&
\Bracket{\wedge^4 \\ \wD^2\\ \wD^0}(-4)\ar[l]_-{\partial_2} & 
0. \ar[l]
}
\]
Our goal is to write the differential $\partial_1$ explicitly.

Let $\{\alpha_1, \dots, \alpha_4\}$ be a basis for $A$, $\{u_{1}, u_{2}\}$ be a
basis for $B_1^*$, and $\{v_1,v_2\}$ be a basis for $B_2^*$. Also, let
$\{v^*_1, v^*_2\}$ be the dual basis for $B_2$. To represent
$\partial_1$ by a matrix, we choose the natural bases of $F(\phi,w)_1$
and $F(\phi,w)_0$ induced by our choice of bases for $A$, $B_1^*$, and
$B_2^*$.  Namely, our basis of $F(\phi,w)_1$ is given by the six
elements of the form
\[
f_{\{i_1,i_2\},\{1,2\}, \emptyset}:=\left( \alpha_{i_1}\wedge
  \alpha_{i_2}\right) \otimes \left( u_1\wedge u_2 \right) \otimes 1,
\]
where $1\leq i_1<i_2\leq 4$.  Our basis of $F(\phi,w)_0$ is given by the
three elements of the form
\[
g_{\emptyset, \emptyset, (j_1,j_2)}:=1\otimes 1\otimes
\left({v_1^*}^{j_1}{v_2^*}^{j_2}\right)
\]
where $(j_1,j_2)\in \mathbb N^2$ and $j_1+j_2=2$.  With notation as in this
subsection, we have that
\begin{equation}\label{eqn:gamma}
  \gamma(f_{\{i_1,i_2\},\{1,2\}, \emptyset} \otimes g_{\emptyset, \emptyset,
    (j_1,j_2)}^*)=\left( \alpha_{i_1}\wedge \alpha_{i_2}\right) \otimes \left( u_1\wedge
    u_2 \right)\otimes \left(v_1^{(j_1)}v_2^{(j_2)} \right).
\end{equation}
If we represent $\phi^{\flat}$ by the matrix of linear forms
\[
\phi^{\flat}=\bordermatrix{
& \alpha_1^* & \alpha_2^* & \alpha_3^* & \alpha_4^* \cr
u_1^*\otimes v_1^* & x_{1, (1,1)} & x_{2,(1,1)} & x_{3,(1,1)} & x_{4,(1,1)} \cr
u_1^*\otimes v_2^* & x_{1, (1,2)} & x_{2,(1,2)} & x_{3,(1,2)} & x_{4,(1,2)} \cr
u_2^*\otimes v_1^* & x_{1, (2,1)} & x_{2,(2,1)} & x_{3,(2,1)} & x_{4,(2,1)} \cr
u_2^*\otimes v_2^* & x_{1, (2,2)} & x_{2,(2,2)} & x_{3,(2,2)} & x_{4,(2,2)} \cr
},
\]
then combining \eqref{eqn:gamma} and \eqref{eqn:inclusionMinors}
allows us to write the image of $\partial_1$ in terms of $2\times 2$
minors of $\phi^{\flat}$.

For example, let us consider the entry of $\partial_1$ corresponding
to $f_{\{1,2\},\{1,2\},\emptyset}$ and $g_{\emptyset, \emptyset,
(2,0)}$. From Appendix~\ref{sec:charfree} we see that
the inclusion
\[
\bigwedge^{2} A \otimes \det B_1^* \otimes \Di^{2}(B_2^*) \subseteq
\bigwedge^{2} A \otimes \bigwedge^{2}(B_1^* \otimes B_2^*)
\]
is given by
\begin{equation}\label{eqn:inclusionMinors}
\begin{cases}
  (u_1 \wedge u_2)\otimes v_1^{(2)} & \mapsto (u_1\otimes v_1)\wedge
  (u_2\otimes v_1)\\ 
  (u_1\wedge u_2)\otimes v_1v_2 & \mapsto (u_1\otimes v_1)\wedge
  (u_2\otimes v_2)+(u_1\otimes v_2)\wedge (u_2\otimes v_1)\\ 
  (u_1\wedge u_2)\otimes v_2^{(2)} & \mapsto (u_1\otimes v_2)\wedge
  (u_2\otimes v_2).
\end{cases}
\end{equation}
Combining \eqref{eqn:gamma} and \eqref{eqn:inclusionMinors}, we
conclude that the entry of $\partial_1$ corresponding to
$f_{\{1,2\},\{1,2\},\emptyset}$ and $g_{\emptyset, \emptyset, (2,0)}$
is given by $(\alpha_{1}\wedge \alpha_{2})\otimes (u_1\otimes v_1)\wedge (u_2
\otimes v_1)$.  Thus, we may write this entry of $\partial_1$ as the
$2\times 2$-minor of $\phi^{\flat}$ obtained by taking the determinant
of the submatrix
\[
\bordermatrix{ & \alpha_1 & \alpha_2 \cr u_1\otimes v_1& x_{1, (1,1)} &
  x_{2,(1,1)} \cr u_2 \otimes v_1 & x_{1, (2,1)} & x_{2,(2,1)} \cr }.
\]
The other entries for $\partial_1$ may be obtained similarly.  See
Example~\ref{ex:422} for a matrix representation of both $\partial_1$
and $\partial_2$ in this example.
\end{example}

\section{Tensor complexes from pinching weights}
\label{sec:main}

We now introduce the notion of pinching weights for a tensor, which enables
us to produce  tensor complexes $F(\phi,w)_\bullet$ that satisfy the
properties of Theorem~\ref{thm:tensor:complexes}.
In contrast with the case of balanced
tensor complexes, there are often many possible pinching weights for a
given tensor $\phi^{a\times \bb}$.

The motivation behind the definition of a pinching weight is the
following.  Recall from \eqref{eqn:cohom:in:each:deg} that the terms
of $F(\phi,w)_\bullet$ can be written as direct sums of certain
cohomology groups on $\PPb$.  Further, since the support $Y(\phi)$ of
$M(\phi,w)$ is independent of $w$
(Corollary~\ref{cor:Mphiw}\eqref{enum:SupportMphiw}), the length of
$F(\phi,w)_\bullet$ is at least $\codim Y(\phi)$, and thus
$F(\phi,w)_\bullet$ must be built from at least this many different
nonzero cohomology groups.  The weight $w$ is a pinching weight
precisely when $F(\phi,w)_\bullet$ is composed of this minimal number
of cohomology groups.

\begin{defn}\label{def:pinching:weight}
A weight vector $w = (w_0, \ldots, w_n) \in\ZZ^{n+1}$ is a \defi{pinching
weight} for $\phi^{a\times \bb}$ if $w_1
<w_2<\cdots<w_n$ and if, for all $1\leq i\leq n$, the intervals
$[w_i+1, w_i+b_i-1]$ lie in $[0,a]$ and are pairwise disjoint.
\end{defn}

The stipulation that $w_1<\cdots<w_n$ is a matter of convention; it can be
guaranteed by permuting the $B_j$.  In addition, we note that if
$F(\phi^{a\times \bb},w)_\bullet$ is a balanced tensor complex, then $w$ is
a pinching weight for $\phi^{a\times \bb}$.

\begin{notation}\label{notation:dw}
  Let $w$ be a pinching weight for $\phi^{a\times \bb}$.  Let $p= a -
  \sum_{j=1}^n (b_j-1)$.  We define two degree sequences and some
  constants in terms of $w$ and the size of $\phi$:
\begin{align*}
  d'(w) & := \left([0,a] \setminus
    \bigcup_{i=1}^n [w_i+1,w_i+b_i-1]\right)\in \ZZ^{p+1}; \\
  d(w) & := d'(w) + (w_0,\dots,w_0)\in \ZZ^{p+1}; \\
  r_i & := \min \{j \mid j > 0 \;\text{and}\; w_j \geq d'_i\}\quad\text{for
    all}\; 0 \leq i \leq p. \qedhere
\end{align*}
\end{notation}

\begin{thm} 
\label{thm:main:structure}
If $w$ is a pinching weight for $\phi$, 
then $F(\phi,w)_\bullet$ is a free complex of length
$a-\sum_{j=1}^n (b_j-1)$, 
and the $i$th term of $F(\phi,w)_\bullet$ is 
\begin{align}
\label{eq:tensor:complex:term}
F(\phi,w)_i &= 
S(-d_i)\otimes \bigwedge^{d_i'} A \otimes 
\bigotimes_{j=1}^{r_i-1} \wD^{d_i'-w_j-b_j} (B_j^*) \otimes 
\bigotimes_{j=r_i}^n \Sc^{w_j-d_i'}(B_j).
\end{align}
The tensor complex $F(\phi,w)_\bullet$ satisfies
Theorem~\ref{thm:tensor:complexes}.  The choice of $G$-equivariant
differential is unique, up to sign.
\end{thm}

\begin{proof}
From the definition of pinching weights and
an argument similar to the proof of 
Proposition~\ref{prop:balanced}, we conclude that
\[
F(\phi,w)_i = S(-d_i) \otimes \bigwedge^{d_i-w_0} A \otimes
\bigotimes_{j=1}^{r_i-1} \rH^{b_j-1}(\PP(B_j),\cO(w_j-d'_i)) \otimes
\bigotimes_{j=r_i}^n \rH^0(\PP(B_j),\cO(w_j - d'_i)).
\] 
This yields \eqref{eq:tensor:complex:term}. 
The desired assertions of
Theorem~\ref{thm:tensor:complexes} then follow from 
minor variants of the arguments in the proofs of 
Propositions~\ref{prop:balanced} and~\ref{prop:canonical},
where we replace the expression \eqref{eqn:balancedMods} by the
expression \eqref{eq:tensor:complex:term}.
The uniqueness (up to sign) of a $G$-equivariant differential follows
by a similar variant of the proof of Proposition~\ref{prop:canonical}.
\end{proof}

\begin{example}\label{ex:722}
Let $a\times \bb=7\times(2,2)$ and $w=(w_0,1,4)$ for any $w_0$.  This 
is a pinching weight for $\phi^{7\times (2,2)}$,
since the intervals
$[w_1+1,w_1+2-1]=[2,2]$ and $[w_2+1,w_2+2-1]=[5,5]$ are disjoint and both
belong to the interval $[0,7]$. 
  The corresponding complex
  $F(\phi,w)_\bullet$ equals the tensor product the complex
  of~\eqref{eqn:threestrand} with $S(-w_0)$.

When $w=(0,-1,6)$, the complex $F(\phi,w)_\bullet$ equals the linear complex
\begin{equation*}
\small \xymatrix @-0.5pc {
\Bracket{\wedge^1\\\wD^0\\ \Sc^5}(-1)
&
\Bracket{\wedge^2\\\wD^1\\ \Sc^4}(-2)\ar[l]
&
\Bracket{\wedge^3\\ \wD^2\\ \Sc^3}(-3)\ar[l]
&
\Bracket{\wedge^4\\ \wD^3\\ \Sc^2}(-4) \ar[l]
&
\Bracket{\wedge^5\\ \wD^4\\  \Sc^1}(-5)\ar[l]
&
\Bracket{\wedge^6\\\wD^5\\ \Sc^0}(-6)\ar[l]
&
0\ar[l]
}.
\end{equation*}

When $w=(-4,1,2)$, the complex $F(\phi,w)_\bullet$ equals:
\begin{equation*}
\small \xymatrix @-0.5pc {
\Bracket{\wedge^0\\ \Sc^1\\ \Sc^2}(4)
&
\Bracket{\wedge^1\\ \Sc^0 \\ \Sc^1}(3)\ar[l]
&
\Bracket{\wedge^4\\\wD^1\\ \wD^0}(0)\ar[l]
&
\Bracket{\wedge^5\\ \wD^2\\ \wD^1}(-1) \ar[l]
&
\Bracket{\wedge^6\\ \wD^3\\ \wD^2}(-2)\ar[l]
&
\Bracket{\wedge^7\\ \wD^4\\ \wD^3}(-3)\ar[l]
&
0\ar[l]
}. \qedhere
\end{equation*}
\end{example}

\begin{rmk}
\label{remark:pinchWtNonCM}
Instead of setting $w$ to be a pinching weight for $\phi$, consider the case
where the intervals $[w_i+1, w_i+b_i-1]$ are pairwise disjoint, but
where we drop the requirement that all the intervals $[w_i+1,
w_i+b_i-1]$ lie in $[0,a]$.  In this case, $M(\phi,w)$ is a
non-Cohen--Macaulay module with a pure resolution. For instance, with
$w=(-4,1,-3)$, $F(\phi^{7\times (2,2)},w)_\bullet$ is
\begin{equation*}
\small \xymatrix @-1pc {
\Bracket{\wedge^0\\ \Sc^1\\ \wD^1}(4)
&
\Bracket{\wedge^1\\ \Sc^0 \\ \wD^2}(3)\ar[l]
&
\Bracket{\wedge^3\\\wD^0\\ \wD^4}(1)\ar[l]
&
\Bracket{\wedge^4\\ \wD^1\\ \wD^5}(0) \ar[l]
&
\Bracket{\wedge^5\\ \wD^2\\ \wD^6}(-1)\ar[l]
&
\Bracket{\wedge^6\\ \wD^3\\ \wD^7}(-2)\ar[l]
&
\Bracket{\wedge^7\\ \wD^4\\ \wD^8}(-3)\ar[l]
& 0 \ar[l]
}.
\end{equation*}
Since $\codim M(\phi,w) = 5$, it is not Cohen--Macaulay.
\end{rmk}

\section{Strands of the Koszul complex}
\label{sec:strands}

We now provide a more elementary description of
$F(\phi,w)_\bullet$ as a complex constructed 
by splicing strands of a Koszul complex together, 
extending the study of matrix complexes in~\cite{buchs-eis} and~\cite[\S A2.6]{eisenbud}. 
The purpose of this section is expository, so
we focus on the example of the universal $7\times (2,2)$ tensor with
pinching weight $w=(0,1,4)$ described in~\eqref{eqn:threestrand}; 
the general case can be treated in a similar fashion.
By Proposition~\ref{thm:main:structure}, 
\[
\beta(F(\phi,w)_\bullet)=\begin{pmatrix}
10&28&-&-&-&-\\
-&-&70&70&-&-\\
-&-&-&-&28&10\\
\end{pmatrix}. 
\]
We now express $F(\phi,w)_\bullet$ in terms of three linear
strands arising from a Koszul complex.  As discussed in
\S\ref{subsec:tensorCxConst}, these are:
\smallskip
\begin{align*}
&\text{ Strand 1:} &\quad & \text{ Strand 2:} &\quad & \text{ Strand 3:}\\
&\xymatrix{ 
\Bracket{\wedge^0\\ \Sc^1\\ \Sc^4} 
&\Bracket{\wedge^1\\ \Sc^0 \\ \Sc^3}(-1) \ar[l] } 
,&& \xymatrix{
\Bracket{\wedge^3\\ \wD^0 \\ \Sc^1}(-3) &\Bracket{\wedge^4\\ \wD^1 \\
\Sc^0}(-4)  \ar[l] } 
,&&
\xymatrix{ \Bracket{\wedge^6\\ \wD^3 \\ \wD^0}(-6) &\Bracket{\wedge^7\\ \wD^4 \\
\wD^1}(-7). \ar[l] } 
\end{align*}
To obtain them from a Koszul complex, we consider the space
$\mathbb A^{a\times \bb}\times \PP(B_1)\times \PP(B_2)$ and let $T:=S\otimes_{\ZZ}
\ZZ[B_1]\otimes_{\ZZ} \ZZ[B_2]$ with the induced $\ZZ^3$-grading.  
If $\mathbf{k} := (k,k,k)\in \ZZ^3$, then  
\[
K(\phi)_\bullet : \xymatrix{ \bigwedge^0 T^7& \bigwedge^1
T^7(-\mathbf{1}) \ar[l] & \dots \ar[l] & \bigwedge^7 T^7(-\mathbf{7})\ar[l]
& 0\ar[l]}
\]
is the $\ZZ^3$-graded Koszul complex from Remark~\ref{rmk:Koszul:forms}. 
For any $\alpha, \beta\in \ZZ$, the subcomplex
$[K(\phi)_\bullet]_{(*,\alpha,\beta)} := \bigoplus_{\gamma \in\ZZ}
[K(\phi)_\bullet]_{(\gamma,\alpha,\beta)}$ of $K(\phi)_\bullet$ is a graded
complex of $S$-modules.  In particular, Strand~1 arises as the $(*,1,4)$
subcomplex of $K(\phi)_\bullet$: 
\begin{align*}
[K(\phi)_\bullet]_{(*,1,4)} 
=& 
\xymatrix{[\bigwedge^0 T^7]_{(*,1,4)}
&[\bigwedge^1 T^7(-\mathbf{1})]_{(*,1,4)}\ar[l]
& \dots \ar[l]
& [\bigwedge^7 T^7(-\mathbf{7})]_{(*,1,4)} \ar[l] & 0 \ar[l]}\\
=&
\xymatrix{
\ \ \ \Bracket{\wedge^0\\ \Sc^1\\ \Sc^4}
\ \ \ \ \ 
&
\ \ \ \ \Bracket{\wedge^1\\ \Sc^0 \\ \Sc^3}(-1)
\ \ \ \ \ 
\ar[l]&0\ar[l]&\dots \ar[l]&0\ar[l]}.
\end{align*}

Strand $2$ also arises from the Koszul complex $K(\phi)_\bullet$, but in a
more subtle manner. 
Let $\wD^\bullet(B_1^*):=\bigoplus_{i=0}^\infty \wD^i(B_1^*)$, which is 
naturally isomorphic as a graded module to the top local
cohomology group of the $\ZZ$-algebra $\ZZ[B_1]$ 
with support in the prime
ideal generated by $B_1$.
Now let $T^{\{1\}}$ be the $T$-module 
$S\otimes \wD^\bullet(B_1^*)\otimes \Sc^\bullet(B_2)$ 
and set $K^{\{1\}}_\bullet:=K(\phi)_\bullet \otimes_T T^{\{1\}}$. 
We then obtain Strand~2 as the $(*,1,4)$ subcomplex of $K^{\{1\}}$:
\small
\begin{align*}
[K^{\{1\}}_\bullet]_{(*,1,4)}&=
\xymatrix{ 
0 
&0 \ar[l] &0 \ar[l] &
\ \ \Bracket{\wedge^3\\ \wD^0\\ \Sc^1}(-3) \ar[l]
\ 
&
\Bracket{\wedge^4\\ \wD^1\\ \Sc^0}(-4) \ar[l]
\ 
\ar[l] &0 \ar[l] &0 \ar[l] 
&0\ar[l]}.
\end{align*}
\normalsize

Finally, Strand~3 is obtained through a similar process. 
If $T^{\{1,2\}}:=S\otimes \wD^\bullet(B_1^*)\otimes \wD^\bullet(B_2^*)$ and
$K^{\{1,2\}}_\bullet:=K(\phi)_\bullet \otimes_T T^{\{1,2\}}$, then Strand~3
arises as the $(*,1,4)$ subcomplex of $K^{\{1,2\}}_\bullet$.

\begin{rmk}\label{rmk:contrasting:strands}
The construction outlined in this section provides a slightly different
view from~\cite[\S A2.6]{eisenbud} 
of building matrix complexes from strands of the Koszul complex,
and we now contrast these approaches.  
Let us consider the case of a $7\times 2$ matrix with
$w=(0,3)$.
The complex $F(\phi,w)_\bullet$ then corresponds to the complex
$\EuScript{C}^2$ of~\cite[\S A2.6]{eisenbud}.  
Incorporating the appropriate twists by determinants into $\EuScript{C}^2$
(as suggested by the footnotes in~\cite[\S A2.6]{eisenbud}), we see
that the complexes $\EuScript{C}^2$ and $F(\phi,w)_\bullet$ are equal. In
both cases, the Betti diagram of the free resolution is
\[
\begin{pmatrix}
3&14&21&-&-&-&-\\
-&-&-&35&42&21&4
\end{pmatrix}.
\]

However, the construction of $F(\phi,w)_\bullet$ differs from the construction
of $\EuScript{C}^2$.  We obtain the first strand of each construction 
in the same manner, as the $(*,2)$ subcomplex of $K(\phi)_\bullet$.  However, the
second strands come from slightly different sources.  

Strand 2 of $\EuScript{C}^2$ is obtained by peeling off the
$(*,3)$ subcomplex of $K(\phi)_\bullet$, which has 
\[
\beta\left( [K(\phi)_\bullet]_{(*,3)}\right) =\begin{pmatrix}
-&-&-&-&-&-&-\\
4&21&42&35&-&-&-
\end{pmatrix},
\]
and then dualizing that strand (and twisting by the appropriate determinants). 
Note that strand $2$ of $\EuScript{C}^2$ originates in homological
degrees $0,1,2$ and $3$ of the complex $K(\phi)_\bullet$, and then duality
is used to turn this strand around.  

By contrast, Strand 2 of
$F(\phi,w)_\bullet$ comes from homological degrees $4,5,6$ and $7$ of a
different complex $K^{\{1\}}_\bullet$:
\[
\beta\left( [K^{\{1\}}_\bullet]_{(*,2)}\right) =\begin{pmatrix}
-&-&-&-&-&-&-\\
-&-&-&35&42&21&4
\end{pmatrix}.
\]
These strands coincide, at least up to a twist by determinants, because of the
self-duality properties of the Koszul complex $K(\phi)_\bullet$.
\end{rmk}

\section{Functoriality properties of tensor complexes}
\label{sec:functoriality}

We now prove Proposition~\ref{prop:functor}, which
describes the functorial properties of the
construction of tensor complexes. We also consider the relation to the
complexes considered in~\cite{BEKS}.

\begin{proof}[Proof of Proposition~\ref{prop:functor}]
  We have $a' \le a$, $w,w'\in \ZZ^{n+1}$, and an inclusion $i \colon
  \ZZ^{a'} \to \ZZ^a$.  First assume that $w=w'$.  This induces a map
  of rings $S' \to S$ (where $S = \ZZ[X^{a\times \bb}]$, $S' =
  \ZZ[X^{a'\times\bb}]$) and a commutative diagram
  \[
  \xymatrix{ \bA^{a\times \bb}\times \PPb\ar[d]_-{\pi}
    \ar[r]^-{\nu}&\bA^{a' \times \bb}\times  \PPb\ar[d]_-{\pi'}\\
    \bA^{a\times \bb}\ar[r]^-{\rho}& \bA^{a' \times \bb}. }
  \]
On $\bA^{a\times \bb}\times \PPb$ and $\bA^{a'\times \bb}\times \PPb$, we have the Koszul complexes $\KK(\phi)_\bullet$ 
and $\KK(\phi')_\bullet$, respectively. The inclusion $i$ induces a natural map 
$\nu^*\KK(\phi')_\bullet \to \KK(\phi)_\bullet$. We thus obtain a natural
  map ${\bf R}\pi_* (\nu^* \KK(\phi')_\bullet) \to {\bf
    R}\pi_*\KK(\phi)_\bullet$.  By the projection formula
  \cite[Proposition II.5.6]{residues}, there is a quasi-isomorphism
   \begin{equation}\label{eqn:qis}
  {\bf R}\pi_* (\nu^* \KK(\phi')_\bullet) \cong 
        \rho^*({\bf R}\pi'_* \KK(\phi')_\bullet)
  \end{equation}
  (noting that ${\bf L}\rho^*$ and ${\bf L}\nu^*$ coincide with
  $\rho^*$ and $\nu^*$, since we apply them to a complex of locally
  free sheaves).
  
  In fact, this map is an isomorphism of complexes. This follows from the
  claim that if $P_\bullet$ and $P'_\bullet $ are minimal (i.e., 
  $\partial P_i \subseteq \fm P_{i-1}$, where $\mathfrak m\subseteq S$
  is the ideal generated by the variables) bounded-below complexes of
  $S$-modules, then a quasi-isomorphism of $P_\bullet$ and
  $P'_\bullet$ induces an isomorphism of these complexes.  To prove
  the claim, we first observe that there is a minimal bounded-below
  complex $\hat P_\bullet$ of free $S$-modules together with maps
  $P_\bullet \leftarrow \hat P_\bullet \rightarrow P'_\bullet$ that realizes the
  quasi-isomorphism. A map between bounded-below projective complexes
  which is a quasi-isomorphism is a homotopy equivalence, and a
  homotopy equivalence between minimal complexes of $S$-modules is an
  isomorphism, proving the claim.

  We thus get a map $\rho^*({\bf R} \pi'_* \KK(\phi')_\bullet) \to
  {\bf R}\pi_* \KK(\phi)_\bullet$. Note that $F(\phi',w)_{\bullet}$ is
  a minimal free resolution in the quasi-isomorphism class of ${\bf R}
  \pi'_* \KK(\phi')_\bullet$, and $F(\phi,w)_{\bullet}$ is a minimal
  free resolution in the quasi-isomorphism class of ${\bf R}\pi_*
  \KK(\phi)_\bullet$.  The above map thus induces the desired map $f_w
  \colon F(\phi',w)_{\bullet} \otimes_{S'} S \to F(\phi,w)_{\bullet}$.

  When $w\ne w'$, we fix a nonzero polynomial $h$ of multidegree
  $w-w'$ on $S'\otimes \ZZ[B_1]\otimes \dots \otimes
  \ZZ[B_n]$, assuming that one exists. Multiplication by $h$
  gives a morphism $\KK(\phi')_{\bullet}(w') \to \KK(\phi')_{\bullet}(w)$.
By taking the global sections of the derived pushforward, we get a morphism
  $F(\phi',w')_{\bullet} \to F(\phi',w)_{\bullet}$. Tensoring with $S$ and
  composing with the map $f_w$ then yields the desired morphism
  $F(\phi',w')_{\bullet} \otimes_{S'} S\to F(\phi,w)_{\bullet}$.
\end{proof}

\begin{example}\label{ex:Func22}
Fix $\bb=(2,2)$ and $w=(0,0,2)$ and consider the tensor complexes
$F(\phi^{a\times \bb}, w)_\bullet$ for $a=2,3$ and $4$.  The Betti diagrams
of these three complexes are
\[
\begin{pmatrix}
3&-&-\\
-&1&-\\
-&-&-
\end{pmatrix},
\begin{pmatrix}
3&-&-\\
-&3&-\\
-&-&-
\end{pmatrix},
\text{ and }
\begin{pmatrix}
3&-&-\\
-&6&-\\
-&-&3
\end{pmatrix},
\]
respectively.  If we choose the multihomogeneous form $1$, then the maps
induced by Proposition~\ref{prop:functor} yield natural inclusions
among these complexes.
\end{example}

Even when $w\ne w'$, the maps induced by Proposition~\ref{prop:functor} can
be simple in special cases.

\begin{lemma}[Homomorphism Push-Forward Lemma]\label{lemma:pushfwdmap}
  With notation as in Proposition~\ref{prop:functor}, assume that
  $F(\phi,w)_\bullet$ is a pure resolution of type $d=(d_0,\dots,
  d_p)$ and that $F(\phi',w')_\bullet$ is a pure resolution of type
  $d'=(d_0', \dots, d_q')$.  Let $h$ be the multihomogeneous form
  determining the morphism of complexes $h_\bullet \colon \mathcal
  K(\phi')_\bullet \to \mathcal K(\phi)_\bullet$ that induces the
  morphism of complexes $\nu_\bullet \colon F(\phi',w')_\bullet \to
  F(\phi,w)_\bullet$.

  Assume further that $d_i=d_i'$ for some $i$, and let  $N:=\sum_{j\leq
  r_i} (b_j-1)$.  Then the map $\nu_i$ may be chosen to be the
  induced map on cohomology, $\rH^N(h_i)$, as in the following diagram:
  \[
  \xymatrix{
    F(\phi',w')_i\ar[rr]^-{\nu_i}\ar[d]_\cong & & F(\phi,w)_i\ar[d]_\cong\\
    \rH^{N}(\mathbb A^{a\times \bb}\times \PPb, \mathcal
    K(\phi')_i)\ar[rr]^-{\rH^N(h_i)}& & \rH^{N}(\mathbb A^{a\times
      \bb}\times \PPb,\mathcal K(\phi)_i). }
  \]
\end{lemma}

\begin{proof}
  As in~\cite[Proof of Proposition~5.3]{ES-JAMS}, we may use the spectral
  sequence with $E^{k,-l}_{1}={\bf R}^{k}\pi_* \mathcal
  K(\phi)_l$ to compute the complex $F(\phi,w)_\bullet$, along with
  a similar spectral sequence to compute
  $F(\phi',w')_\bullet$.  We thus construct $\nu_\bullet$
  by considering the map induced by $h_\bullet$ on these
  spectral sequences.  Since $d_i=d_i'$, one may check that on the
  $E_1$ page, the induced map in position $(N,-i)$ is given by
  \[
  \rH^{N}(h_i) \colon \rH^{N}(\mathbb A^{a\times \bb}\times \PPb,
  \mathcal K(\phi')_i)\to \rH^{N}(\mathbb A^{a\times \bb}\times
  \PPb,\mathcal K(\phi)_i).
  \]
  Since the terms of the complexes $F(\phi,w)_\bullet$ and
  $F(\phi',w')_\bullet$ come from the $E_1$ page of this spectral sequence,
  this map may be chosen as the map $\nu_i \colon
  F(\phi',w')_i\to F(\phi,w)_i$.
\end{proof}

\begin{rmk}\label{rmk:Func22}
The proof of~\cite[Theorem~1.2]{BEKS} can be reinterpreted in terms of
Proposition~\ref{prop:functor} and Lemma~\ref{lemma:pushfwdmap}.  Namely,
the pure resolutions of~\cite[Theorem~1.2]{BEKS} are
specializations of certain tensor complexes of the form $a\times
(2,2,\dots,2)$, and the morphisms constructed between two such resolutions are
of the form of the morphisms given by
Proposition~\ref{prop:functor}.  
However, we note that Proposition~\ref{prop:functor} 
does not directly imply that result, 
as the above map of complexes could be null-homotopic.  
The essential step in the proof of~\cite[Theorem~1.2]{BEKS}
is checking that certain maps of complexes induce
nonzero maps $M(\phi^{a'\times\bb},w')\to M(\phi^{a\times\bb},w)$, which
requires analyzing the detailed description of $\nu_\bullet$ provided by
Lemma~\ref{lemma:pushfwdmap}.
\end{rmk}

\section{Properties of the module $M(\phi,w)$}
\label{sec:Mphi}
The goal of this section is to prove Corollary~1.5 and
Propositions~\ref{prop:resultant} and \ref{prop:indecomp}.
We begin by discussing some facts about the support
$Y(\phi)$ of $M(\phi,w)$.  In \S\ref{sec:MDVs} we explore the geometry of
$Y(\phi)$ further.

Recall the diagram of~\eqref{eqn:Yphi}.  The scheme $Y(\phi)$ is
integral since it is the scheme-theoretic image of the integral scheme
$Z(\phi)$.  Throughout this section we identify $\bA^{a\times \bb}$
with the space of $\ZZ$-linear maps $\psi \colon \bB^* \to A^*$. 
For a linear subspace $V$ of $\bB^*$ we write
$[V]$ for the corresponding subspace in $\PP(\bB)$. So for any map
$\psi\in \mathbb A^{a\times \bb}$, we may think of $[\ker(\psi)]$ as a
linear subspace of $\PP(\bB)$.  Let $\Seg(\bB)$ denote the image of
the Segre embedding $\PPb \to \PP(\bB)$.

\begin{prop}\label{prop:YphiKerCapSeg}
  The annihilator of $M(\phi,w)$ is the prime
  ideal that defines the integral scheme $Y(\phi)$. 
  Under the identification $\PPb \cong \Seg(\bB)$, we have
  \[
Z(\phi) = \{(\psi, y ) \in \bA^{a \times \bb} \times \PPb \mid
y \in [\ker(\psi)]\}.
\]
We therefore have
\[
Y(\phi)=\{\psi\in \Hom(\bB^*, A^*) \mid [\ker(\psi)]\cap \Seg(\bB) \neq
\varnothing\}\subseteq \mathbb A^{a\times \bb}.
\]
\end{prop}
\begin{proof}
  The first assertion follows from ~\cite[Theorems~5.1.2(b), 5.1.3(a)]{weyman}, which imply that $M(\phi,w)$ is a module over the
  normalization of $Y(\phi)$. Since $Z(\phi)$ is the total space of
  $\cS = \sheafHom((\bB^* \otimes \cO_{\PPb}) /
  \cO_{\PPb}(-\mathbf{1}), A^* \otimes \cO_{\PPb})$, we may think of
  $Z(\phi)$ as the set of maps $\psi \colon \bB^*\otimes \cO_{\PPb} \to
  A^*\otimes \cO_{\PPb}$ whose kernel contains a rank 1 tensor, yielding
  the second assertion. The final assertion is now immediate.
\end{proof}

\begin{rmk}\label{rmk:resultant}
We now explain how Proposition~\ref{prop:YphiKerCapSeg} implies 
Proposition~\ref{prop:resultant}, which states that $Y(\phi)$ may be
interpreted 
as a resultant variety for multilinear equations on $\PPb$. 
As in Remark~\ref{rmk:Koszul:forms}, we view a point in 
$\mathbb A^{a\times \bb}$ as a collection 
$\widetilde{\mathbf{f}}=(\widetilde{f}_1, \dots, \widetilde{f}_a)$ 
of multilinear forms on $\PPb$. 
Then Proposition~\ref{prop:YphiKerCapSeg} implies that $Z(\phi)$ 
is the incidence variety
\[
Z(\phi)=\{ (\widetilde{\mathbf{f}},y)\in \bA^{a \times \bb} \times \PPb \mid
y \in V_{\PPb}(\widetilde{f}_1, \dots, \widetilde{f}_a)\}, 
\]
and thus it follows that $Y(\phi)$ has the resultant interpretation 
\[
Y(\phi)=\{ \widetilde{\mathbf{f}}\in \bA^{a \times \bb} \mid
 V_{\PPb}(\widetilde{f}_1, \dots, \widetilde{f}_a)\ne \emptyset\}.
\]
This yields Proposition~\ref{prop:resultant} because, 
for any pinching weight $w$, the complex $F(\phi,w)_\bullet$ 
resolves a module whose support equals $Y(\phi)$.  Hence 
the minors of $\partial_1$ cut out $Y(\phi)$ set-theoretically 
by~\cite[Proposition 20.7]{eisenbud}.
\end{rmk}

\begin{rmk} \label{rmk:birational} The map $\mu \colon Z(\phi) \to
  Y(\phi)$ restricts to an isomorphism over the (possibly empty) open
  subset of $Y(\phi)$ consisting of those $\psi$ such that $\ker \psi$
  contains a rank 1 tensor that is unique up to scalar multiple. Now,
  if $a > \sum_{i=1}^n (b_i - 1)$ (i.e., if $\codim Y(\phi) \ge 1$ by
  Corollary~\ref{cor:Mphiw}\eqref{enum:SupportMphiw}), then such maps
  $\psi$ exist, and hence $\mu$ is a birational morphism. We see this
  as follows. If $b_1 \cdots b_n \le a$, then we may choose any rank 1
  tensor and define $\psi$ to be a map whose kernel is spanned by the
  chosen rank 1 tensor. Now suppose that $b_1 \cdots b_n > a$.  Then,
  since $a > \sum_i (b_i - 1)$, there is a $(b_1\cdots
  b_n-1-a)$-plane (i.e., a linear subvariety of codimension $a$) in
  $\PP(\bB)$ that intersects $\Seg(\bB)$ in exactly one
  point. Define $\psi$ to be a map with this linear subvariety as its
  kernel.
\end{rmk} 
 
We are now prepared to prove Proposition~\ref{prop:indecomp} and Corollary~\ref{cor:Mphiw}. 

\begin{proof}[Proof of Proposition~\ref{prop:indecomp}]
Note first that the sheaf 
\[
\widetilde{M}(\phi,w)= \mu_* 
\left(\cO_{Z(\phi)} \otimes \pi_2^*\left(\cO_{\PPb}(w_1, \dots,
w_n)\right)\right)\otimes \cO(-w_0)
\]
is a twist of the pushforward of a line bundle.  Using the fact
that $\mu$
is birational by Remark~\ref{rmk:birational}, 
we see that  $\widetilde{M}(\phi,w)$ is a rank-one sheaf on
$Y(\phi)$.  Since $Y(\phi)$ is irreducible and
$M(\phi,w)$ is Cohen--Macaulay, the indecomposability of $M(\phi,w)$ follows
immediately.
\end{proof}

\begin{proof}[Proof of Corollary~\ref{cor:Mphiw}]
For~\eqref{enum:SupportMphiw}, the fact that the support of
$M(\phi,w)$ does not depend on $w$ follows from
Proposition~\ref{prop:YphiKerCapSeg}. The codimension formula follows
from~\eqref{eq:tensor:complex:term} and~\cite[Theorem~5.1.6(a)]{weyman} (while this result is proven when the base ring is a field of characteristic 0, we may reduce to this case because 
$Y(\phi)$ is flat over $\ZZ$). 

For~\eqref{enum:MphiwGenPerf}, recall that $M(\phi,w)$ is 
Cohen--Macaulay and has a uniformly minimal resolution over $\ZZ$
(Theorem~\ref{thm:tensor:complexes}\eqref{item:tensor:CM} 
and \eqref{item:tensor:flat}).
By uniform minimality, $\Tor^{\ZZ}_1(M(\phi,w), \ZZ/\ell) = 0$ and $\ell M(\phi,w) \neq M(\phi,w)$
for all primes $\ell$. In other words, $M(\phi,w)$ is a faithfully flat
$\ZZ$-module (see, for example, \cite[Theorems~7.2,~7.8]{matsumura}), and hence $M(\phi,w)$ is generically perfect.

For~\eqref{enum:MphiwMult}, we may assume that $w_0=0$.  Set $p :=
a - \sum_j (b_j-1)$, so $d(w) = d'(w) = (d_0', \dots, d_p')$.  To simplify the notation,
we use $M$ to denote $M(\phi,w)$ throughout the rest of this proof.  Since
$M$ is a Cohen--Macaulay
module with a pure resolution of type $d(w)$
(by Theorem~\ref{thm:tensor:complexes}\eqref{item:tensor:CM} 
and~\eqref{item:tensor:pure}), we see
from~\cite[Theorem~1.2]{huneke-miller} that
\[
e(M) = 
\frac{1}{(\codim M)!}\left( \prod_{i=1}^n (d_i'-d_0')\right) 
\beta_{0,d_0'}(M).
\]
(Huneke and Miller prove \cite[Theorem~1.2]{huneke-miller} only for cyclic
Cohen--Macaulay modules, but \cite[(1.3)]{huneke-miller} can be modified
by multiplying by $\beta_{0,d_0'}(M)$ to make the proof work for all
Cohen--Macaulay modules with pure resolutions.)

We use \eqref{eq:tensor:complex:term} to compute
$\beta_{0,d'_0}(M)$.  Recall that $r_0 = \min \{ j \mid w_j \ge
d'_0 \}$ and 
observe that, since $w$ is a pinching weight, we have $[d_0',a]=\{d_0',
d_1', \dots, d_n'\}\sqcup  \bigsqcup_{j\geq r_0} [w_j+1, w_j+b_j-1]$.
We may rewrite this as 
$[0,a-d_0']=\{0,d_1'-d_0', \dots, d_n'-d_0'\}\sqcup  \bigsqcup_{j\geq r_0} [w_j+1-d_0', w_j+b_j-1-d_0'].$ Since 
\[
\rank_{\ZZ} \left( \bigotimes_{j=r_0}^n \Sc^{w_j-d_0'}(B_j)\right) = \prod_{j \ge r_0}
\frac{(w_j+1-d'_0) \cdots (w_j + b_j - 1 - d'_0)}{(b_j-1)!},
\]
this yields 
\begin{align*}
 \rank_{\ZZ} \left(\bigotimes_{j=r_0}^n \Sc^{w_j-d_0'}(B_j) \right)
  \cdot \prod_{j=1}^n (d_j'-d_0') & = \frac{(a-d_0')!}{\prod_{j\geq r_0} (b_j-1)!}.
\end{align*}
In addition, we have $[0,d_0'-1]=\bigsqcup_{j<r_0} [w_j+1,w_j+b_j-1]$, 
since $w$ is a pinching weight.  Multiplying by $-1$ and adding
$d_0'$, we obtain the equality $[1,d_0']=\bigsqcup_{j<r_0}
[d_0'-w_j-b_j+1,d_0'-w_j-1], $ and we similarly see that
\[
\rank_{\ZZ}\left(\bigotimes_{j=1}^{r_0-1} \wD^{d_0'-w_j-b_j} (B_j^*)
\right)=\frac{d_0!}{\prod_{j<r_0} (b_j-1)!}.
\]
Finally, we combine these to get the multiplicity of $M$:
\begin{align*}
  e(M)&=
  \frac{1}{(\codim M)!}\left( \prod_{i=1}^n (d_i'-d_0')\right)\cdot 
  \left(\rank_{\ZZ} F(\phi,w)_0\right)\\
  &=\frac{1}{(\codim M)!}\binom{a}{d_0'}
  \left(\frac{d_0!}{\prod_{j<r_0} (b_j-1)!} \right)
  \left(\frac{(a-d_0')!}{\prod_{j\geq r_0} (b_j-1)!} \right)\\
  &=\frac{1}{(\codim M)!}\frac{a!}{\prod_{j=1}^n (b_j-1)!}.\qedhere
\end{align*}
\end{proof}
See Remark~\ref{rmk:orderES} for a surprising consequence of the above
formula for $e(M(\phi,w))$.

\begin{rmk}\label{rmk:self-dual}
By imposing symmetry, we can obtain tensor complexes that are equivariantly self-dual.  For example, reconsider the tensor complex from \eqref{eqn:threestrand}.  Based on the representations that arise in the free resolution, the complex exhibits certain symmetries; but it is not a self-dual resolution of $S$-modules.  

However, a variant of this complex is self-dual.  Let $\Bbbk=\ZZ[\frac{1}{2}]$.  Since $B_1\cong B_2$, we may identify these free modules and consider $\Sc^2(B_1)\otimes_{\ZZ} \Bbbk\subseteq B_1\otimes B_2\otimes_{\ZZ} \Bbbk $.  Let $S':=\Bbbk[A\otimes \Sc^2(B_1)\otimes_{\ZZ} \Bbbk]$ and $\phi'$ be the universal symmetric tensor in $A\otimes \Sc^2(B_1)\otimes S'$.  By applying the above inclusion, we may view $\phi'$ as a tensor in $A \otimes B_1\otimes B_2\otimes S'$ and thus construct $F(\phi',w)_\bullet$ as a complex of $S'$-modules.  

The complex $F(\phi',w)_\bullet$ is equivariantly self-dual as a complex of $S'$-modules. This self-duality is forced by the uniqueness of equivariant differentials, as discussed in \S\ref{sec:differentials}.  A similar construction works whenever $B_i\cong B_{n-i}$ for all $i$ and $w_j+w_{n+1-j}=-b_j$ for all $j$.
\end{rmk}

\section{Hyperdeterminantal varieties}\label{sec:MDVs}
There are two special cases where the supporting variety 
$Y(\phi)$ has been previously studied in some detail.
First, if there is a unique $i$ such that $b_i\ne 1$, then
$Y(\phi)$ is the determinantal variety defined by the
maximal minors of a universal matrix.  
Motivated by this example, we
refer to $Y(\phi)$ as a \defi{hyperdeterminantal variety}.
The second case where hyperdeterminantal varieties have previously
been studied is when $\codim Y(\phi)=1$.  
As we prove in Proposition~\ref{prop:hyperdeterminant}, in this case, 
 $Y(\phi)$ is defined by a hyperdeterminant of the boundary format.

 The main goal of this section is to prove Theorem~\ref{thm:hyperdeterminants2}, as well as describe other geometric properties
 of hyperdeterminantal varieties.  Based on the
two special cases above, one might wonder if the variety $Y(\phi)$ is
 Cohen--Macaulay in general. This turns out to be entirely false:
 Proposition~\ref{prop:MDV:TFAE} shows that $Y(\phi)$ is
 Cohen--Macaulay if and only if it is either a determinantal variety,
 a hypersurface, or all of $\mathbb A^{a\times \bb}$.  We 
 consider the singular locus of $Y(\phi)$ in 
 Proposition~\ref{prop:normal:singular}; in the hyperdeterminantal case, our 
 result recovers a portion of~\cite[Theorem~0.5(a)]{wz}.

To begin with hyperdeterminantal hypersurfaces, 
the tensor $\phi^{a\times \bb}$ is said to have the \defi{boundary format} 
when $a -\sum_{i=1}^n (b_i-1)=1$~\cite[\S 14.3]{gkz}. 
In this case, there is a corresponding hyperdeterminant
$\Delta_{a\times \bb}$, which is generally defined over a field of
characteristic $0$.  However, since $\Delta_{a\times \bb}$ 
is unique up to scalar multiple, we view it as a
polynomial over $\ZZ$ that is not divisible by any prime number $\ell$, 
so that it is unique up to sign.

\begin{prop}\label{prop:hyperdeterminant}
  Let $\phi = \phi^{a\times \bb}$ be of the boundary format and $w$ be
  any pinching weight for $\phi$. Then
  $F(\phi,w)_\bullet$ is a free resolution of length $1$, and hence
  $\partial_1$ is a square matrix.  Up to sign, the hyperdeterminant $\Delta_{a\times \bb}$
  equals $\det (\partial_1)$.
\end{prop}
\begin{proof}
  We first show that $Y(\phi)$ equals the vanishing of the
  hyperdeterminant $\Delta_{a\times \bb}$.  By
  Corollary~\ref{cor:Mphiw}, we may choose any $w$ to compute
  $Y(\phi)$.  We set $w_0:=0$, $w_1:=1$, and $w_i:=( \sum_{j<i} b_i) -
  (i-2)$ for $i\geq 1$.  We confirm that this yields a pinching weight
  for $\phi$ by computing
\[
[w_1+1, w_1+b_1-1]=\begin{cases}
[2,b_1] & i=1,\\
[(3-i)+\sum_{j<i} b_j,(5-i)+\sum_{j\leq i} b_j] & i\geq 2.
\end{cases}
\]
By Theorem~\ref{thm:main:structure}, the resulting free resolution is a two
term linear complex:
\[
\small \xymatrix @-0.5pc { \Bracket{\wedge^0 \\ \Sc^{w_1} \\ \Sc^{w_2}
    \\ \vdots \\ \Sc^{w_n}} && \Bracket{\wedge^1 \\ \Sc^{w_1-1} \\
    \Sc^{w_2-1} \\ \vdots \\ \Sc^{w_n-1}}(-1)\ar[ll]_-{\partial_1} &
  0\ar[l]},
\]
where $\partial_1$ is a $G$-equivariant map.  

The source and target of $\partial_1$ can naturally be associated with
the source and target of the matrix $\partial_A$ from
\cite[Proposition~14.3.2]{gkz}, which is used to compute the
hyperdeterminant $\Delta_{a\times \bb}$.  Clearly $\partial_A$ is
$G$-equivariant by definition.  We claim that $\partial_1$ and
$\partial_A$ differ by $\pm 1$.  After passing to $\QQ$, we see (by an
argument similar to Lemma~\ref{lemma:multfree}) that the map of
representations $[\partial_1]_1\colon [F_1]_1\otimes \QQ\to
[F_0]_1\otimes \QQ$ is an injective map from an irreducible
representation to a multiplicity-free representation.  A similar
statement holds for $[\partial_A]_1$, and hence $[\partial_1]_1$ and
$[\partial_A]_1$ differ by an integer scalar.  Hence it follows that
$\det(\partial_1)$ is an integral scalar multiple of $\Delta_{a\times
  \bb}$.  However, since $Y(\phi)$ is irreducible, it follows that
$\det(\partial_1)$ is also, up to sign, a power of an irreducible
polynomial.  This proves that $\det(\partial_1)$ and $\Delta_{a\times
  \bb}$ are equal, up to sign.

Now let $w$ be any pinching weight for $\phi$, and let $\partial_1$ be
the corresponding differential on the $2$-term complex.  Since
$Y(\phi)$ does not depend on $w$, $\det(\partial_1)$
is a power of $\Delta_{a\times \bb}$.  Since $M(\phi,w)$ is
Cohen--Macaulay of codimension $1$, its multiplicity equals the 
degree of $\det(\partial_1)$.  By combining
Corollary~\ref{cor:Mphiw}\eqref{enum:MphiwMult} and
\cite[Corollary~14.2.6]{gkz}, it follows that $\deg
\det(\partial_1)=\deg \Delta_{a\times \bb}$, completing the proof.
\end{proof}

We note that \cite[Theorem 14.3.1]{gkz} provides a resultant interpretation 
for a hyperdeterminant of the boundary format. 
As discussed in Remark~\ref{rmk:resultant}, 
this interpretation generalizes to higher codimension, 
enabling us to prove Theorem~\ref{thm:hyperdeterminants2}.

\begin{proof}[Proof of Theorem~\ref{thm:hyperdeterminants2}]
As it is enough to show this result after passing 
to an algebraically closed field $\Bbbk$, we replace $Y(\phi)$, etc., 
by their corresponding objects over $\Spec(\Bbbk)$. 
By Remark~\ref{rmk:resultant}, we may then apply the resultant interpretation 
of $Y(\phi)$ to view the $\Bbbk$-points of $Y(\phi)$ 
as systems of multilinear equations $\widetilde{\mathbf{f}}$ that have 
a nonempty vanishing locus in $\PPb$.  

Recall that $a'=1+\sum_{i=1}^n (b_i-1)$, 
and let $I$ be the ideal of $a'\times \bb$ hyperdeterminants 
from~\eqref{eq:hyperdet ideal}.  We claim that set-theoretically, $V(I)=Y(\phi)$. 
Note that $Y(\phi)\subseteq V(I)$, since any collection of $a'$ polynomials 
in the vector space $\langle \widetilde{f}_1, \dots, \widetilde{f}_a\rangle$ 
must have a common root, and thus all of the corresponding hyperdeterminants must vanish by~\cite[Theorem~14.3.1]{gkz}.

For the reverse inclusion, suppose that there exists a point 
$\widetilde{\mathbf{f}} \in V(I)\setminus Y(\phi)$. 
We thus have that $\widetilde{\mathbf{f}}$ has no common zero in $\PPb$. 
Since $V(I)$ and $Y(\phi)$ are both $G$-equivariant, 
we may assume after a $\GL(A)$-change of coordinates that 
$\widetilde{f}_1, \dots, \widetilde{f}_{a'-1}$ intersect 
in a finite number of points $\{P_1, \dots, P_t\}\in \PPb$.  
We now consider the vector space 
$W:=\langle \widetilde{f}_{a'}, \dots, \widetilde{f}_a\rangle$ 
and choose $\widetilde{g}\in W$. 
Since every hyperdeterminant of every sub-tensor $\phi'$ of $\phi$ of size $a'\times \bb$ 
vanishes on $\widetilde{\mathbf{f}}$, 
there must be some $P_i$ that is a root of $\widetilde{g}$. 
Consequently, the incidence locus 
\[
\{ (\widetilde{g},P_i) \in W\times \{P_1, \dots, P_t\} \mid \widetilde{g}(P_i)=0\} 
\]
is a closed sublocus of $W\times \{P_1, \dots, P_t\}$ that surjects onto $W$. 
It then follows that there is some connected component 
of this incidence locus that alone surjects onto $W$; 
in other words, there is some $P_i$ that is simultaneously a root 
of all polynomials in $W$. 
This $P_i$ is then also a common zero of $\widetilde{\mathbf{f}}$, contradicting our assumption that $\widetilde{\mathbf{f}}\notin Y(\phi)$.
\end{proof}

\begin{rmk}\label{rmk:sturmfels}
Bernd~Sturmfels has pointed out that $Y(\phi^{a\times b\times 2})$ has a
second interpretation as a resultant variety as well.
For simplicity, we work over a field $\Bbbk$.
By identifying points of $\bA^{a\times b\times 2}_{\Bbbk}$ with maps in
$\Hom(\Bbbk^a \otimes \Bbbk^2,\Bbbk^{b})$, we may think of a point $\psi\in
\bA^{a \times b\times 2}_{\Bbbk}$ as a linear map
\[
\iota_{\psi}\colon \PP^{b-1}\to \PP^{2a-1}.
\]
The image of $\iota_{\psi}$ then intersects the Segre variety $\PP^{a-1}\times \PP^1$ if and only if $\psi$ belongs to $Y(\phi)$.
This can be checked directly as follows.  Let $U_1, \dots, U_{b}$ be a
sequence of $2\times a$ matrices which span the image of $\iota_{\psi}$.  The
image of $\iota_{\psi}$ intersects the Segre variety if and only if there
exist nontrivial scalars $\lambda_i$ and $(\alpha_1,\alpha_2)$ such that
$(\alpha_1,\alpha_2)$ belongs to the kernel of $\sum_{i=1}^{b}
\lambda_iU_i$.  This is equivalent to the statement that the rank $1$
tensor $(\lambda_i\alpha_j)\in \Bbbk^{b}\otimes \Bbbk^{2}$ belongs to the
kernel of $\psi^{\flat}$, which is equivalent to $\psi\in Y(\phi)$ by
Proposition~\ref{prop:YphiKerCapSeg}.
\end{rmk}

We now provide a more detailed description of the geometry of
$Y(\phi)$. When $b_i > 1$ for only one index $i$,
$Y(\phi)$ is a determinantal variety defined by the maximal minors of
a matrix of indeterminates. We thus investigate the situation when $b_i > 1$ for at
least two indices $i$.

\begin{prop} 
\label{prop:MDV:TFAE}
Suppose that $b_i > 1$ for at least two indices $i$ and that $Y(\phi) \ne
\bA^{a \times \bb}$. Then $Y(\phi)$ is not normal. If additionally 
$\codim Y(\phi) \geq 2$, then $Y(\phi)$ is not Cohen--Macaulay.
\end{prop}

\begin{proof} 
  Since $\mu \colon Z(\phi) \to Y(\phi)$ is birational
  (Remark~\ref{rmk:birational}), it suffices, by Zariski's
  connectedness theorem, to show that there is a fiber of $\mu$ that
  is not geometrically connected.  
  
  Let $\psi \in Y(\phi)$ be a generic map.  We claim that $\ker(\psi)
  \cap \Seg(\bB)$ is a single point $x$.  If $b_1 \cdots b_n \le a$,
  then $\ker \psi$ is 1-dimensional; therefore, the intersection is a
  single point.  If $b_1 \cdots b_n > a$, then the kernel of a map
  $\psi \colon \bB^* \to A^*$ has codimension $a$.  Since $a >
  \dim \Seg(\bB)$ and $[\ker \psi] \cap Y(\phi) \neq \varnothing$, we
  obtain the claim.
  
  Let $\Bbbk$ be the algebraic closure of the residue field of $\psi$,
  so that $x$ is $\Bbbk$-rational. Pick an additional $\Bbbk$-rational
  point $y$ on $\Seg(\bB)$ but not on $[\ker \psi]$ such that the line
  joining $x$ and $y$ does not lie in $\Seg(\bB)$. (Here we use
  the hypothesis that $b_i > 1$ for at least two $i$. Note that if
  $b_i > 1$ for at most one $i$, then $\Seg(\bB)$ is a linear
  subvariety of $\PP(\bB)$.)  Pick a basis for $\bB^*$ containing
  $x$ and $y$, and let $\psi'$ be a map that agrees with $\psi$ on all basis
  elements except $y$ and sends $y$ to $0$. Then $\psi' \in Y(\phi)$ and
  $[\ker \psi']$ intersects $\Seg(\bB)$ in finitely many points (but at
  least two). Hence the fiber over $\psi'$ is not geometrically
  connected.

  Now assume that $\codim Y(\phi) = a - \sum_{i=1}^n (b_i-1)\geq
  2$. Then, by Proposition~\ref{prop:normal:singular}, $Y(\phi)$ is regular
  in codimension one.  By the Serre criterion for
  normality~\cite[Theorem 11.5]{eisenbud}, $Y(\phi)$ does not satisfy
  the condition ($S_2$), so is not Cohen--Macaulay.
\end{proof}

The following proposition provides a multilinear analogue of the
classical fact that the singular locus of a determinantal variety
consists of those maps whose kernel has dimension higher than the
generic value.

\begin{prop}
  \label{prop:normal:singular} 
  Suppose that $b_i > 1$ for at least two indices $i$ and that
  $Y(\phi) \ne \bA^{a \times \bb}$.  Then the singular locus
  $Y(\phi)_{\rm{sing}}$ of $Y(\phi)$ coincides with the non-normal
  locus $Y(\phi)_{\rm nn}$ of $Y(\phi)$.  In particular,
  \[
  Y(\phi)_{\rm{sing}}= \{\psi\in Y(\phi) \mid [\ker(\psi)] \cap
  \Seg(\bB) \; \text{is not a single reduced point}\}.
  \]
  Furthermore, $Y(\phi)_{\rm{sing}}$ is irreducible of codimension
  $a-\sum_{i=1}^n(b_i-1)$ in $Y(\phi)$.  
\end{prop}

\begin{proof}
  Let $Y_1 \defeq \{\psi\in Y(\phi) \mid [\ker(\psi)] \cap \Seg(\bB)
  \;\text{is not a single reduced point}\}$. We first show that $Y_1$
  is irreducible.  Let
  $\Delta \subset \Seg(\bB) \times \Seg(\bB)$ be the diagonal
  subscheme and $U := (\Seg(\bB) \times \Seg(\bB)) \minus \Delta$.
  Write $q_1$ and $q_2$ for the two projection morphisms $\Seg(\bB)
  \times \Seg(\bB) \to \Seg(\bB)$. Note that $\LL := (q_1^*\cO(-1,
  \dots, -1) \oplus q_2^*\cO(-1, \dots, -1))|_U$ is naturally a
  subbundle of the trivial bundle $\bB^* \otimes \cO_U$.  Let $Z'$ be
  the total space of $\sheafHom((\bB^* \otimes \cO_U)/ \LL, A^*\otimes
  \cO_U)$; note that $Z'$ is an irreducible subvariety of $\bA^{a
    \times \bb} \times U$, which is the total space of
  $\sheafHom(\bB^* \otimes \cO_U, A^*\otimes \cO_U)$.
A point $\psi \in \bA^{a
    \times \bb}$ lies in the image of  $Z'$ if and only if $[\ker \psi]\cap \Seg(\bB)$ consists of more than one point.

Hence, every point of $Y_1$ lies in the closure of the image of $Z'$ (which is irreducible), except possibly the loci
of $\psi$ such that $[\ker(\psi)]\cap \Seg(\bB)$ consists of a single nonreduced point.  Thus, to complete our argument
that $Y_1$ is irreducible, we must show that any such $\psi$ lies in the closure of $Y_1$.
Fix some $\psi_0$ such that $[\ker \psi_0] \cap \Seg(\bB)$ is a single nonreduced point,
and write $[\ker \psi_0]$ as a sum of lines $L_1 + L_2 +
  \cdots + L_r$ such that $L = L_1$ is a tangent line to $\Seg(\bB)$
  at $x$. Since a tangent line at a smooth point is a limit of secant
  lines, there is a family of secant lines $L_t$ that have $L$ as
  their limit, and we write $H_t := L_t + L_2 + \cdots + L_r$.  There
  is then a compatible family of $\psi_t$ such that $[\ker \psi_t] = H_t$
  and $\psi_t$ limits to $\psi_0$.  Since $L_t$ is a secant line, it follows that $H_t\cap \Seg(\bB)$
  is supported on more than point, and hence $\psi_t\in Y_1$.  Since $\psi_0$ is in the closure of the family $\psi_t$,
  it follows that $\psi_0$ also lies in $Y_1$, as desired.
 
We next compute the codimension of $Y_1$ in $Y(\phi)$. The map $Z' \to Y_1$ is a 2-to-1 map over the dense
  open subset of $Y_1$ where $[\ker \psi]$ intersects
  $\Seg(\bB)$ in two points.  Therefore
  \begin{align*}
    \dim Y_1& = \dim Z' \\
    & = 2\dim \Seg(\bB) + a(b_1 b_2 \cdots b_n -2) \\
    & = \dim Y(\phi) - \rank A + \dim \Seg(\bB).
  \end{align*}
  Hence the codimension of $Y_1$ in $Y(\phi)$ is $a -
  \sum_{i=1}^n (b_i - 1)$.

  Finally, we claim that $Y_1$ coincides with both the singular locus
  $Y(\phi)_{\rm{sing}}$ and the non-normal locus $Y(\phi)_{\rm{nn}}$
  of $Y(\phi)$.  As noted in Remark~\ref{rmk:birational}, $\mu \colon
  Z(\phi)\to Y(\phi)$ is birational over the open set
  $Y(\phi)\setminus Y_1$.  Since $Z(\phi)$ is smooth, we see that
  $Y(\phi)_{\rm{sing}} \subseteq Y_1$.  Now, since $b_i > 1$ for at
  least two indices $i$, $\Seg(\bB)$ is not a linear subvariety of
  $\PP(\bB)$. Thus, as argued in the proof of
  Proposition~\ref{prop:MDV:TFAE}, there exist $\psi \in Y_1$ such
  that $[\ker(\psi)] \cap \Seg(\bB)$ set-theoretically consists of at
  least two reduced points.  Since $Y_1$ is irreducible, it follows
  that a general point of $Y_1$ has this property.  Any such point is
  a non-normal point of $Y(\phi)$, and since both $Y_1$ and
  $Y(\phi)_{\rm{nn}}$ are closed, we conclude that $Y_1 \subseteq
  Y(\phi)_{\rm{nn}}$.  Of course, the non-normal locus always sits in
  the singular locus, and we thus obtain the chain
  \[
  Y(\phi)_{\rm{sing}} \subseteq Y_1\subseteq Y(\phi)_{\rm{nn}}\subseteq Y(\phi)_{\rm{sing}},
  \]
proving that these loci coincide.
\end{proof}

\begin{rmk} In the case when $\phi^{a \times \bb}$ is a tensor of the 
  boundary format, Proposition~\ref{prop:normal:singular} recovers the first
  part of \cite[Theorem~0.5(a)]{wz}, which says that the singular
  locus of a hyperdeterminantal hypersurface (of the boundary format) 
  is irreducible and has codimension 1. 
  In these cases, $Y(\phi)$ is Cohen--Macaulay
  since it is a hypersurface, but it fails to be normal.
\end{rmk}

\begin{rmk}\label{rmk:CMmodules}
  A conjecture of M.~Hochster asserts that every complete local domain
  has a finitely generated maximal Cohen--Macaulay module~\cite[Conjecture~6, p.10]{HochCBMS75}. This is known to be true in only
  a handful of cases~\cite{HochBigCMWitt75, GrifSplittingCM98,
    KatzSmallCM99, ScheDimFiltr99}.  By combining
  Theorem~\ref{thm:tensor:complexes} and
  Proposition~\ref{prop:MDV:TFAE}, we can construct finitely generated
  maximal Cohen--Macaulay modules $M(\phi, w)$ with
  non-Cohen--Macaulay supports $Y(\phi)$.  At all points $y$ where the
  completion of $\cO_{Y,y}$ is a domain (i.e., at the unibranched
  points of $Y$) we get new examples where Hochster's conjecture
  holds.  As far as we know, these examples are not covered by any
  previously known results. For instance, we could take $y$ to be the
  $\ZZ/p$-point lying over the origin of $\bA^{a\times \bb}$.
\end{rmk}

\begin{example}\label{ex:rank}
Consider the case $a\times \bb=3\times (2,2)$ and $w=(0,0,1)$.  
Then $F(\phi,w)_\bullet$ is a two-term complex
$S^2(-3)\overset{\partial_1}{\to} S^2$.
By the method for writing out $\partial_1$ described in \S\ref{subsec:flattenings}, 
we see that each entry of $\partial_1$ 
corresponds to a specific $3\times 3$ minor of $\phi^{\flat}$.  

Now, let $\widetilde{\phi}\in \CC^3\otimes \CC^2\otimes \CC^2$ denote
a $\CC$-point of $\mathbb A^{a\times \bb}$.  By
\cite[Theorem~1.1]{landsberg-weyman}, the border rank of the tensor
$\widetilde{\phi}$ is less than $3$ if and only if the $3\times 3$
minors of $\phi^{\flat}$ vanish when evaluated at $\widetilde{\phi}$.
This is equivalent to asking that the specialization of $\partial_1$
at $\widetilde{\phi}$ yields the zero matrix.  Thus in this case, the
border rank of the tensor $\widetilde{\phi}$ is determined by the
homological properties of the specialization of the tensor complex.
It would be interesting to study whether similar connections hold in
more generality.
\end{example}

\section{Eisenbud--Schreyer pure resolutions are balanced tensor complexes}
\label{sec:ESpures}

The existence of pure resolutions of type $d$ for an arbitrary degree
sequence $d$ was originally conjectured
in~\cite[Conjecture~2.4]{boij-sod1}.  The first construction of such 
pure resolutions in arbitrary
characteristic appears in~\cite[\S5]{ES-JAMS}.
Theorem~\ref{thm:ESpures} below implies that each of these
Eisenbud--Schreyer pure resolutions can be realized as the
specialization of some balanced tensor complex.  Each
of these resolutions is constructed from a 
sequence of sufficiently generic multilinear forms 
$\bg := g_1, \dots, g_a$ on $\bA^n_{\Bbbk}\times \PPb$, where $\Bbbk$
is any field; set $R:=\Bbbk[x_1, \dots, x_n]$ and 
denote the corresponding pure resolution of $R$-modules by $\ES(\bg,d)_\bullet$.

\begin{thm}\label{thm:ESpures}
  Let $d = (d_0, \ldots, d_n)$ be a degree sequence, and
  $\ES(\bg,d)_\bullet$ be an Eisenbud--Schreyer pure resolution.  Let
  $a:=d_n-d_0$, $b_i:=d_i-d_{i-1}$, and $w:=(d_0,0,d_1,d_2,\dots,
  d_{n-1})$.  Then there exists a map $\ZZ[X^{a\times \bb}] \to R$
  such that
  \[
  \ES(\bg,d)_\bullet\cong F(\phi^{a\times\bb},w)_\bullet
  \otimes_{\ZZ[X^{a \times {\bf b}}]} R.
  \]
\end{thm}

\begin{proof}
Since each $g_i$ is multilinear, we may write $g_i=\sum_{J} g_{i,J}y_J$,
where the $g_{i,J}$ are linear forms on $\bA^n$ and where $y_J$ is a multilinear form
on $\PPb$.  We then define a map $\ZZ[X^{a\times \bb}]\to R$ by 
$x_{i,J}\mapsto g_{i,J}.$  This yields a commutative diagram:
\[
\xymatrix{
  \bA^n\times \PPb\ar[d]_-{\pi'}\ar[r]^-{\nu}&\bA^{a\times \bb}\times
  \PPb\ar[d]_-{\pi}\\ 
  \bA^n\ar[r]^-{\rho}& \bA^{a\times \bb}. }
\]
By the projection formula \cite[Proposition II.5.6]{residues}, we
get a quasi-isomorphism
\[ 
{\bf R}\pi'_*(\nu^* \KK(\phi)_\bullet) \cong \rho^*({\bf R}\pi_*
\KK(\phi)_\bullet)
\]
(noting that ${\bf L}\rho^*$ and ${\bf L}\nu^*$ coincide with $\rho^*$
and $\nu^*$, since we apply them to a complex of locally free
sheaves). The argument immediately following \eqref{eqn:qis} yields an
isomorphism of complexes. Using the notation of
Remark~\ref{rmk:Koszul:forms}, we have $\nu^*(f_i)=g_i$, so $\nu^*
\KK(\phi)_\bullet$ is the Koszul complex used
in~\cite[Theorem~5.1]{ES-JAMS} to construct the complex
$\ES(\bg,d)_\bullet$.
\end{proof}

\begin{rmk}\label{rmk:ESoverZZ}
In~\cite[Proposition~5.2]{ES-JAMS}, Eisenbud and Schreyer illustrate explicit multilinear forms over $\ZZ$ that satisfy the necessary
genericity conditions.  We note that the Theorem~\label{thm:ESpures} also holds when $R=\ZZ[x_1, \dots, x_n]$ and, in this case, $\ES(\bg,d)_\bullet$ is a uniformly minimal resolution of a generically perfect module $M$ of codimension $n$.
\end{rmk}

\begin{rmk}\label{rmk:orderES}
By combining Corollary~\ref{cor:Mphiw}\eqref{enum:MphiwMult} and
Theorem~\ref{thm:ESpures}, we recover the curious fact that the
multiplicity of the Eisenbud--Schreyer pure resolution of type
$d=(d_0,\dots, d_n)$ depends only on the unordered(!) set of first
differences $\{d_1-d_0, \dots, d_n-d_{n-1}\}$.  We first learned of this
fact through a conversation with Eisenbud and Schreyer.
\end{rmk}

\section{New Families of Pure Resolutions}
\label{sec:new:pures}

We have shown that a tensor $\phi^{a\times \bb}$ 
and a pinching weight $w$ yield a pure resolution 
$F(\phi^{a\times \bb},w)_\bullet$ of type $d(w)$ 
(Notation~\ref{notation:dw}). 
Informally, we may think of this as a map $(a,\bb,w)\mapsto d(w)$, 
where $w$ is a pinching weight for $\phi^{a\times\bb}$. 
From this perspective, the proof of Theorem~\ref{thm:family}
describes the fibers of this map.

\begin{proof}[Proof of Theorem~\ref{thm:family}]
  Let $d \in \ZZ^{p+1}$.  We will describe all the choices of $a,\bb$,
  and pinching weight $w$ such that $F(\phi^{a\times \bb},w)_\bullet$
  is a pure resolution of type $d$.  (The module $M(\phi, w)$ is
  Cohen--Macaulay by Theorem~\ref{thm:main:structure}.)  Let $c \leq
  d_0$ and $C \geq d_p$ be integers, and view $d$ as a subsequence of
  $\{c,c+1, \dots, C\}$. Subdivide $\{c,c+1, \dots, C\} \setminus
  \{d_0, \ldots, d_p\}$ into sequences $s^{(j)}$ of
  consecutive integers, where $1 \leq j \leq n$.  
  We may assume that $\min(s^{(j+1)}) > \min(s^{(j)})$ for all $j$.

  Let $a := C - c$ and $b_j := |s^{(j)}| + 1$ for $1 \leq j \leq
  n$. (Here $|\cdot|$ denotes the length of the sequence.) Let $w_0 :=
  c$ and $w_j := \min(s^{(j)}) - c - 1$, $1 \leq j \leq n$.  Since
  $s^{(j)} = \{w_j+c+1,\dots,w_j+b_j+c-1\}$, we see that the intervals
  $[w_j+1, w_j+b_j-1]$ are disjoint and contained in $[0,a]$.
  Therefore $w$ is a pinching weight
  (Definition~\ref{def:pinching:weight}) for $\phi$.  Note that by
  construction, $d(w)=d$. Thus we have chosen $a$, $\bb$, and $w$ so
  that $F(\phi^{a\times\bb},w)_\bullet$ is a pure resolution of type
  $d$, and there are infinitely many such choices.
\end{proof}

\begin{rmk}\label{rmk:EFWconj}
If $w$ is a pinching weight for $\phi$ (so that $F(\phi,w)_\bullet$ is a
pure resolution of type $d(w)$), then the Betti diagram of
$F(\phi,w)_\bullet$ is an integral multiple of the Betti diagram of the
corresponding Eisenbud--Schreyer pure resolution.  In particular,
Theorem~\ref{thm:family} has no implications
for~\cite[Conjecture~6.1]{efw}.
\end{rmk}
\begin{table}[h]
\caption{Pure resolutions of type $d = (0,3)$ with parameters $c = -2$ and $C = 4$.}
\begin{tabular}{|c|c|c|c|}\hline
Subdivision& $a\times \bb$ & $w$&$\beta(F(\phi,w)_\bullet)$ \\ \hline
$(-2,-1), (1,2), (4)$ & $6\times (3,3,2)$ & $(-2,-1,2,5)$ &
$ \begin{pmatrix}
60&-\\
-&-\\
-&60 
\end{pmatrix}$  \\ \hline
$(-2), (-1), (1,2), (4)$ & $6\times (2,2,3,2)$ & $(-2,-1,0,2,5)$ &
$ \begin{pmatrix}
120&-\\
-&-\\
-&120 
\end{pmatrix}$
\\ \hline
$(-2,-1), (1), (2), (4)$ & $6\times (3,2,2,2)$ & $(-2,-1,2,3,5)$ &
$ \begin{pmatrix}
120&-\\
-&-\\
-&120 
\end{pmatrix}$
\\ \hline
$(-2), (-1), (1), (2), (4)$ & $6\times (2,2,2,2,2)$ & $(-2,-1,0,2,3,5)$&
$ \begin{pmatrix}
240&-\\
-&-\\
-&240
\end{pmatrix}$  \\ \hline
\end{tabular}
\label{tab:04}
\end{table}

\begin{example}
\label{example:03}
Consider the degree sequence $d=(0,3)$.  Table~\ref{tab:04}
illustrates the various constructions of pure resolutions of type $d$
with $c=-2$ and $C=4$.
\end{example}

\begin{example}
The complexes $F_\bullet$ and $F'_\bullet$ 
in \cite[Example~6.5]{BEKS} are also specializations of tensor complexes; this follows from an argument similar to the proof of
Corollary~\ref{thm:ESpures}.  
Namely, the complex $F_\bullet$ is a specialization of the tensor complex
for an $8\times (2,2,2,2)$ tensor with $w = (0,0,2,6,7)$; the complex
$F_\bullet'$ is a specialization of the tensor complex for a $7\times
(2,2,2,2)$ tensor with $w' = (0,-1,2,4,5).$ 
We obtain 
\begin{equation*}
  \small \xymatrix @-0.5pc {
    F_\bullet: \quad 
    {\begin{array}{c} \Bracket{\wedge^0\\ \Sc^0\\
        \Sc^2\\ \Sc^6\\ \Sc^7}\end{array}} 
    &
    {\begin{array}{c} \Bracket{\wedge^2\\\wD^0
        \\ \Sc^0 \\ \Sc^4\\ \Sc^5}(-2)\ar[l] \end{array}}
    &
   {\Bracket{\wedge^4\\ \wD^2
        \\ \wD^0\\ \Sc^2\\ \Sc^3}(-4)\ar[l]}
    &
    {\Bracket{\wedge^5\\ \wD^3 \\
        \wD^1\\ \Sc^1\\ \Sc^2}(-5) \ar[l] }
    &
    {\Bracket{\wedge^6\\\wD^4\\
        \wD^2\\ \Sc^0\\ \Sc^1}(-6)\ar[l] } 
     &
     0 \ar[l]
      }
\end{equation*}
and 
\begin{equation*}
\small \xymatrix @-0.5pc {
F'_\bullet: \quad 
\Bracket{\wedge^1\\\wD^0\\ \Sc^1\\ \Sc^3\\ \Sc^4}(-1)
&
\Bracket{\wedge^2\\\wD^1\\ \Sc^0 \\ \Sc^2\\ \Sc^3}(-2)\ar[l]
&
\Bracket{\wedge^4\\ \wD^3\\ \wD^0\\ \Sc^0\\ \Sc^1}(-4)\ar[l]
&
\Bracket{\wedge^7\\ \wD^6\\ \wD^3\\ \wD^1\\ \wD^0}(-7)\ar[l]
&
0. \ar[l]
}
\end{equation*}
The nonzero map between these resolutions is induced by the natural
inclusion $A'\subseteq A$ whose cokernel is the final summand of $\ZZ$
in $A$.  See also Remark~\ref{rmk:Func22}.
\end{example}

\section{Detailed Example of a Tensor Complex}
\label{sec:detailed}

\begin{example}\label{ex:422}
Let $\phi$ be the universal $4\times (2,2)$ tensor, and $w=(0,0,2)$.  
We consider the complex $F(\phi,w)_\bullet$.
This is one of the simplest
examples of a tensor complex which is not a matrix complex. The resulting
complex $F(\phi^{4\times(2,2)},(0,0,2))_\bullet$ is
\[
\xymatrix{
\Bracket{ \wedge^0 \\ \Sc^0\\ \Sc^2}&
\Bracket{\wedge^2 \\ \wD^0\\  \Sc^0}(-2)\ar[l]_{\partial_1}&
\Bracket{\wedge^4 \\ \wD^2\\ \wD^0}(-4)\ar[l]_{\partial_2} & 
0 \ar[l]
},
\]
which has the Betti diagram 
\[
\begin{pmatrix}
3&-&-\\
-&6&-\\
-&-&3
\end{pmatrix}.
\]

To describe the differentials $\partial_1$ and $\partial_2,$ we first write
the flattening $\phi^{\flat} \colon A^*\to B_1^* \otimes B_2^*$:
\[
\phi^{\flat}=\bordermatrix{
& 1 & 2 & 3 & 4 \cr
a & x_{1, (1,1)} & x_{2,(1,1)} & x_{3,(1,1)} & x_{4,(1,1)} \cr
b & x_{1, (1,2)} & x_{2,(1,2)} & x_{3,(1,2)} & x_{4,(1,2)} \cr
c & x_{1, (2,1)} & x_{2,(2,1)} & x_{3,(2,1)} & x_{4,(2,1)} \cr
d & x_{1, (2,2)} & x_{2,(2,2)} & x_{3,(2,2)} & x_{4,(2,2)} \cr
}.
\]
For $I\subseteq \{a,b,c,d\}$ and $J\subseteq \{1,2,3,4\}$ with $|I|=|J|$,
we denote the corresponding minor of $\phi^{\flat}$ by $\Delta_{I;J}$. For
instance, $\Delta_{ab;12}$ is the $2\times 2$ minor from the upper left
corner of $\phi^{\flat}$. 

We set $a_1, \dots, a_4$ as a basis of $A$,
$u_1, u_2$ a basis of $B_1$, and $v_1, v_2$ a basis of $B_2$.  Following
the notation and the method of \S\ref{subsec:flattenings}, we then obtain
\begin{align*}
\partial_1^T& =
\bordermatrix{
& g_{\varnothing, \varnothing, (2,0)}& 
g_{\varnothing, \varnothing, (1,1)} & 
g_{\varnothing, \varnothing, (0,2)}  \cr
f_{\{1,2\},\{1,2\},\varnothing}&\Delta_{ac;12} & \Delta_{ad;12}+\Delta_{bc;12} & \Delta_{bd;12}\cr
f_{\{1,3\},\{1,2\},\varnothing}&\Delta_{ac;13} & \Delta_{ad;13}+\Delta_{bc;13} & \Delta_{bd;13}\cr
f_{\{1,4\},\{1,2\},\varnothing}&\Delta_{ac;14} & \Delta_{ad;14}+\Delta_{bc;14} & \Delta_{bd;14}\cr
f_{\{2,3\},\{1,2\},\varnothing}&\Delta_{ac;23} & \Delta_{ad;23}+\Delta_{bc;23} & \Delta_{bd;23}\cr
f_{\{2,4\},\{1,2\},\varnothing}&\Delta_{ac;24} & \Delta_{ad;24}+\Delta_{bc;24} & \Delta_{bd;24}\cr
f_{\{3,4\},\{1,2\},\varnothing}&\Delta_{ac;34} & \Delta_{ad;34}+\Delta_{bc;34} & \Delta_{bd;34}\cr
} \\
\end{align*}
\vspace{-.7cm}
\begin{align*}
\text{and} \qquad 
\partial_2 & =
\bordermatrix{
& e_{\{1234\},(2,0),\varnothing}& 
e_{\{1234\},(1,1),\varnothing} & 
e_{\{1234\},(0,2),\varnothing}  \cr
f_{\{1,2\},\{1,2\},\varnothing}& \Delta_{ab;34} & (\Delta_{ad;34}-\Delta_{bc;34}) & \Delta_{cd;34}\cr
f_{\{1,3\},\{1,2\},\varnothing}& -\Delta_{ab;24} & -(\Delta_{ad;24}-\Delta_{bc;24}) & -\Delta_{cd;24}\cr
f_{\{1,4\},\{1,2\},\varnothing}& \Delta_{ab;23} & (\Delta_{ad;23}-\Delta_{bc;23}) & \Delta_{cd;23}\cr
f_{\{2,3\},\{1,2\},\varnothing}& \Delta_{ab;14} & (\Delta_{ad;14}-\Delta_{bc;14}) & \Delta_{cd;14}\cr
f_{\{2,4\},\{1,2\},\varnothing}& -\Delta_{ab;13} & -(\Delta_{ad;13}-\Delta_{bc;13}) & -\Delta_{cd;13}\cr
f_{\{3,4\},\{1,2\},\varnothing}& \Delta_{ab;12} & (\Delta_{ad;12}-\Delta_{bc;12}) & \Delta_{cd;12}\cr
}.
\end{align*}

The fact that each entry of $\partial_1\partial_2$ equals zero follows from
a generalized Laplace expansion of a singular matrix.  For instance, let us
consider the $(1,1)$ entry of $\partial_1\partial_2$, which is given by
\[
(\partial_1\partial_2)_{1,1}=\Delta_{ac;12}\Delta_{ab;34} - \Delta_{ac;13} \Delta_{ab;24} 
+ \Delta_{ac;14} \Delta_{ab;23} + \Delta_{ac;23} \Delta_{ab;14} 
- \Delta_{ac;24} \Delta_{ab;13} + \Delta_{ac;34} \Delta_{ab;12}.
\]
By the generalized Laplace expansion
formula~\cite[\S 1.6]{northcott}, this equals the determinant of
\[
\bordermatrix{
& 1 & 2 & 3 & 4 \cr
a & x_{1, (1,1)} & x_{2,(1,1)} & x_{3,(1,1)} & x_{4,(1,1)} \cr
c & x_{1, (2,1)} & x_{2,(2,1)} & x_{3,(2,1)} & x_{4,(2,1)} \cr
a & x_{1, (1,1)} & x_{2,(1,1)} & x_{3,(1,1)} & x_{4,(1,1)} \cr
b & x_{1, (1,2)} & x_{2,(1,2)} & x_{3,(1,2)} & x_{4,(1,2)} \cr
}.
\]
But the above matrix has a repeated row, and hence this determinant is
zero.  Similar arguments show that all entries of $(\partial_1\partial_2)$
equal $0$.
\end{example}

\begin{example}\label{ex:422support}
Continuing with the $4\times 2\times 2$ example above, we 
compute the defining ideal of $Y(\phi)$.
In order to use representation theory and computations from 
{\tt Macaulay2}~\cite{M2}, 
we work over $\QQ$ instead of $\ZZ$.
From the presentation matrix $\partial_1$ for $M(\phi,w)$, 
we compute directly in {\tt Macaulay2} that $Y(\phi)$ 
is defined by $1$ quartic and $10$ sextic equations.  The quartic equation arises as the determinant of $\phi^{\flat}$, which corresponds to the subrepresentation 
\[
\bS_{1,1,1,1}(A)\otimes \bS_{2,2}(B_1^*)\otimes
\bS_{2,2}(B_2^*)\subseteq \Sc^4(A \otimes B_1^* \otimes B_2^*).
\]
The sextic equations correspond to the hyperdeterminants of all 
$3\times 2\times 2$ subtensors of $\phi$ and arise as the subrepresentation 
\[
\bS_{2,2,2}(A)\otimes \bS_{3,3}(B_1^*)\otimes
\bS_{3,3}(B_2^*)\subseteq \Sc^6(A \otimes B_1^* \otimes B_2^*).
\]

These equations all have geometric significance. 
Namely, as discussed in Remark~\ref{rmk:resultant}, 
$Y(\phi)$ parametrizes quadruples of multilinear forms $(f_1, \dots, f_4)$ 
on $\PP^1\times \PP^1$ with $V(f_1, \dots, f_4)\ne \emptyset$. 
Since the $H^0(\PP^1\times \PP^1,\cO(1,1))$ is $4$-dimensional and 
base-point-free, 
the vector space $\langle f_1, \dots, f_4\rangle$ has dimension at most $3$. 
This explains the presence of the quartic $\det(\phi^{\flat})$.

In addition, if $V(f_1, \dots, f_4)\ne \emptyset$, then 
$V(g_1, g_2, g_3)\ne \emptyset$ for every triplet 
$g_1, g_2, g_3\in \langle f_1, \dots, f_4\rangle$. 
For a such a triplet, $V(g_1,g_2,g_3)\ne \emptyset$ if and only if 
its corresponding $3\times 2\times 2$ hyperdeterminant vanishes. 
Applying this to all $3\times 2\times 2$ subtensors yields the 
$10$-dimensional space of sextic equations.
\end{example}

\appendix 
\section{Characteristic-free multilinear algebra}
\label{sec:charfree}
We review some characteristic-free multilinear algebra. See
\cite[\S 1.1]{weyman} and \cite{abw}\footnote{The first formula in
\cite[p.247]{abw} is a multiple of the second formula in {\it loc.\ cit.} 
and does not have desirable characteristic-free properties.} for more 
details.

Let $E$ be a finitely generated $\ZZ$-module and $d$ a
positive integer. Let $\Sigma_d$ denote the symmetric group on $d$
letters. The \defi{symmetric power} $\Sc^d(E)$ is the quotient of
$E^{\otimes d}$ by the submodule generated by elements of the form
$e_1 \otimes \cdots \otimes e_d - e_{\sigma(1)} \otimes \cdots \otimes
e_{\sigma(d)}$ for $\sigma \in \Sigma_d$. The \defi{divided power}
$\Di^d(E)$ is the submodule of $\Sigma_d$-invariants of $E^{\otimes
  d}$. We have a canonical isomorphism $\Di^d(E^*) = \Sc^d(E)^*$. The
\defi{exterior power} $\bigwedge^d E$ is the quotient of $E^{\otimes
  d}$ by the submodule generated by elements of the form $e_1 \otimes
\cdots \otimes e_d - {\rm sgn}(\sigma) e_{\sigma(1)} \otimes \cdots
\otimes e_{\sigma(d)}$ for $\sigma \in \Sigma_d$, where ${\rm
  sgn}(\sigma)$ is the determinant of $\sigma$ when written as a
permutation matrix. One could also define the exterior power as a
submodule of $E^{\otimes d}$, but in accordance with
Remark~\ref{rmk:super}, one must make a distinction between the two
when $E$ is a $\ZZ/2$-graded module. If $E$ is a free $\ZZ$-module,
then each module defined is also a free $\ZZ$-module.

For each of the three definitions above, one can take direct sums
over all $d \ge 0$, and the resulting modules can be given the
structure of a Hopf algebra. In particular, they are equipped with a
multiplication $m$ and comultiplication $\Delta$, which we will make
use of.

Now for $E$ and $F$ free $\ZZ$-modules of finite rank, we
define the following inclusions.
\begin{asparaenum}
\item $\bigwedge^d E \otimes \bigwedge^d F \to \Sc^d(E \otimes F)$ 
is defined by mapping $e_1 \wedge \cdots \wedge e_d \otimes
  f_1 \wedge \cdots \wedge f_d$ to the determinant of the
  matrix $(e_i \otimes f_j)_{i,j=1,\dots,d}$.
\item $\Phi_d \colon \bigwedge^d E \otimes \Di^d F \to \bigwedge^d(E
  \otimes F)$ will be defined by induction on $d$.  For the base case,
  set $\Phi_1$ to be the identity.  For $d>1$, extend linearly the map on elements of the form 
  $x = e_1 \wedge \cdots \wedge e_d \otimes f_1^{(\alpha_1)} \cdots f_r^{(\alpha_r)}$, 
  where $\alpha_1 + \cdots + \alpha_r = d$, given by
  \[
  \Phi_d(x) := \sum_{i=1}^r (e_1 \otimes f_i) \wedge \Phi_{d-1}(e_2
  \wedge \cdots \wedge e_d \otimes f_1^{(\alpha_1)} \cdots
  f_i^{(\alpha_i - 1)} \cdots f_r^{(\alpha_r)}).
  \]
\item $\Di^d E \otimes \Di^d F \to \Di^d(E \otimes F)$ is the dual of 
the map $\Sc^d(E^* \otimes F^*) \to \Sc^d(E^*) \otimes
  \Sc^d(F^*)$, which is given by
  $(e_{i_1} \otimes f_{j_1})
  \cdots (e_{i_d} \otimes f_{j_d}) \mapsto (e_{i_1} \cdots
  e_{i_d}) \otimes (f_{j_1} \cdots f_{j_d})$.
\end{asparaenum}

\section{Schur functors in characteristic zero}
\label{sec:schurfunctors}
We review some representation theory of $G = \GL_n(\QQ)$ and
Schur--Weyl duality.  See \cite[\S 2]{weyman} and \cite[\S 5]{kraft} for general background.

A sequence of nonnegative integers $\lambda = (\lambda_1, \dots,
\lambda_n)$ is a partition if $\lambda_1 \ge \lambda_2 \ge \cdots \ge
\lambda_n$. If $\lambda_n \ne 0$, then $n$ is the length of
$\lambda$. We set $|\lambda| := \lambda_1
+ \cdots + \lambda_n$, and write $\lambda \vdash
|\lambda|$. Write $1^d$ for the partition consisting of $d$ 1's. Given
two partitions $\lambda$ and $\mu$, we write $\lambda \subseteq \mu$
if $\lambda_i \le \mu_i$ for all $i$. If $\lambda\subseteq\mu$, we say
that $\mu / \lambda$ is a \defi{horizontal strip}, if $\mu_1 \ge
\lambda_1 \ge \mu_2 \ge \lambda_2 \ge \cdots \ge \mu_n \ge \lambda_n$,
denoted $\mu / \lambda \in \mathrm{HS}$. Define $\lambda'$ to be
the partition such that $\lambda'_i = \#\{ j \mid \lambda_j \ge i\}$.
Given $\lambda \subseteq \mu$, we say that $\mu / \lambda$ is a
\defi{vertical strip} if $\mu' / \lambda'\in\mathrm{HS}$, denoted 
$\mu / \lambda \in \mathrm{VS}$.
 
Let $E$ be the $n$-dimensional vector representation of $G$.  The
finite-dimensional irreducible polynomial representations of $G$ are
indexed by partitions $\lambda$ of length at most $n$, and a general
finite-dimensional polynomial representation of $G$ is a direct sum of
irreducible representations. Let $\bS_\lambda E$ denote the
irreducible representation corresponding to $\lambda$, using the
convention that ${\bf S}_\lambda E = 0$ if $\lambda_{n+1} > 0$. 
In particular, $\Sc^d E = {\bf S}_{(d)} E$ 
and $\bigwedge^d E = {\bf S}_{1^d} E$. 

Pieri's rule gives tensor product decompositions
\begin{align} \label{eqn:pieri} \bS_\lambda E \otimes \Sc^d E \cong
  \bigoplus_{\substack{\mu \vdash |\lambda| + d\\ \mu / \lambda \in
      \mathrm{HS}}} {\bf S}_\mu E
  \qquad \text{and}\qquad 
  \bS_\lambda E \otimes
  \bigwedge^d E \cong \bigoplus_{\substack{\mu \vdash |\lambda| + d \\
      \mu / \lambda \in \mathrm{VS}}} \bS_\mu E.
\end{align}
See \cite[Corollary 2.3.5]{weyman} (there, $L_\lambda E$ is isomorphic
to our $\bS_{\lambda'}E$). These formulas remain valid if we replace
$E$ by its dual $E^*$.

Let $\Sigma_k$ be the symmetric group on $k$ letters. There are
commuting actions of $G$ and $\Sigma_k$ on $E^{\otimes k}$.
\defi{Schur--Weyl duality} \cite[Proposition 5.9]{kraft} is the $G
\times \Sigma_k$-equivariant decomposition
\[
E^{\otimes k} \cong \bigoplus_{\substack{\lambda \vdash k \\
    \lambda_{n+1} = 0}} \bS_\lambda E \otimes \chi_\lambda,
\]
where $\chi_\lambda$ are irreducible representations of $\Sigma_k$.
We use that $\chi_{(k)}$ is the trivial representation of
$\Sigma_k$, $\chi_{(1^k)}$ is the one-dimensional sign representation, 
and more generally, $\chi_\lambda \otimes \chi_{(1^k)}
= \chi_{\lambda'}$.

We use the following consequence of Schur--Weyl duality. 
Let $E_1, \dots, E_r$ be vector spaces 
and consider $\bS_\lambda(E_1 \otimes \cdots\otimes E_r)$ 
as a representation of $\GL(E_1) \times \cdots \times\GL(E_r)$. 
The irreducible representations of 
$\GL(E_1) \times \cdots\times \GL(E_r)$ 
are indexed by $r$-tuples of partitions, so 
\begin{align} \label{eqn:schurweyl1}
\bS_\lambda(E_1 \otimes \cdots \otimes E_r) \cong \bigoplus_{\mu^1,
  \dots, \mu^r} (\bS_{\mu^1} E_1 \otimes \cdots \otimes {\bf
  S}_{\mu^r} E_r)^{\oplus g_{\lambda, \mu^1, \dots, \mu^r}}
\end{align}
for some nonnegative integers $g_{\lambda, \mu^1, \dots, \mu^r}$ (the
{\bf Kronecker coefficients}).  We now apply Schur--Weyl duality to
$(E_1 \otimes \cdots \otimes E_r)^{\otimes k}$ in two different ways,
where $k = |\lambda|$.  First, we have
\[
(E_1 \otimes \cdots \otimes E_r)^{\otimes k} \cong \bigoplus_{\nu
  \vdash k} \bS_\nu(E_1 \otimes \cdots \otimes E_r) \otimes
\chi_\nu
\]
as $\GL(E_1 \otimes \cdots \otimes E_r) \times\Sigma_k$-representations. 
Second, we have
\[
E_1^{\otimes k} \otimes \cdots \otimes E_r^{\otimes k} \cong \left(
  \bigoplus_{\mu^1 \vdash k} \bS_{\mu^1} E_1 \otimes \chi_{\mu^1}
\right) \otimes \cdots \otimes \left( \bigoplus_{\mu^r \vdash k}
{\bS}_{\mu^r} E_r \otimes \chi_{\mu^r} \right)
\]
as $\GL(E_1) \times \cdots \times \GL(E_r)
\times\Sigma_k$-representations.  Restricting to the action of
$\Sigma_k$ and comparing the $\chi_\lambda$-isotypic component of both
expressions, we see that $g_{\lambda, \mu^1, \dots, \mu^r}$ is the
multiplicity of $\chi_\lambda$ in the product $\chi_{\mu^1}
\otimes \cdots \otimes \chi_{\mu^r}$.  Since all
representations of $\Sigma_k$ are self-dual, this yields
\begin{align} \label{eqn:schurweyl2} 
g_{\lambda, \mu^1, \dots, \mu^r} =
  \dim (\chi_\lambda \otimes \chi_{\mu^1} \otimes \cdots \otimes
  \chi_{\mu^r})^{\Sigma_k},
\end{align}
where the superscript indicates that invariants are taken.  In light
of \eqref{eqn:schurweyl2}, $g_{\lambda, \mu^1, \dots, \mu^r}$ is
invariant under permutation of all of its indices.  In particular, we
deduce the Cauchy identities
\begin{equation} \label{eqn:cauchy}
\begin{split}
  \Sc^d(E_1 \otimes E_2) &\cong \bigoplus_{\lambda \vdash d} {\bf
    S}_\lambda E_1 \otimes \bS_\lambda E_2 
  \qquad \text{and}\qquad 
  \bigwedge^d(E_1
  \otimes E_2) \cong \bigoplus_{\lambda \vdash d} \bS_\lambda E_1
  \otimes \bS_{\lambda'} E_2.
\end{split}
\end{equation}

\end{document}